\documentclass[preprint,12pt]{elsarticle}

\usepackage{hyperref}
\usepackage[ruled,linesnumbered]{algorithm2e}
\usepackage{float}
\usepackage{algpseudocode}
\usepackage{caption}
\usepackage{amssymb}
\usepackage{verbatim} 
\usepackage{amsmath,bm}
\usepackage{graphicx}
\usepackage{natbib}
\usepackage{graphicx}
\usepackage{subcaption}
\usepackage{hyperref}
\usepackage{tabularx}
\usepackage[titletoc,title]{appendix}
\newcommand{\repeatedfootnote}[2]{%
  \hypertarget{fn:#1}{}%
  \hyperlink{fn:#1}{\footnotemark[#1]}%
  \footnotetext[#1]{#2}%
}

\usepackage[ruled,linesnumbered]{algorithm2e}
  
  \DeclareMathOperator*{\argmin1}{arg\,min}
\setcitestyle{authoryear,open={(},close={)}} 

\usepackage{tabularx}
\usepackage{array}
\usepackage{amsmath,bm}
\usepackage{siunitx}
\usepackage{amsmath}
\usepackage{epstopdf}
 \usepackage{graphicx}
 \usepackage{booktabs}
 \usepackage{rotating}
\usepackage{tabularx}
 \usepackage{tikz}
 \usepackage{comment}

 \usepackage{longtable}
 \usepackage{xtab}
\usepackage{mathtools,nccmath}
 \usepackage{nicefrac}
 \usepackage{float}
 \usepackage{footnote}
 \usepackage{xcolor}
  \newcommand{\set}{\mathcal}
\newlength\myindent
\usepackage[titletoc,title]{appendix}
\usepackage{xr}
\usepackage{longtable}
\usepackage{amsmath}
\usepackage{graphicx}
\usepackage{xcolor}

\graphicspath{ {./images/} }
\usepackage[ruled,linesnumbered]{algorithm2e}
\usepackage{framed}
\usepackage{caption}
\usepackage{subcaption}
\usepackage[justification=centering]
{caption}
\journal{}

\begin{document}

\begin{frontmatter}



\title{PAMSO: Parametric Autotuning Multi-time Scale Optimization Algorithm}


\author[inst1]{Asha Ramanujam}

\affiliation[inst1]{organization={Davidson School of Chemical Engineering},
            addressline={Purdue University, 480 W. Stadium Ave, West Lafayette, IN, 47907 USA}}

\author[inst1]{Can Li}

\begin{abstract}
Optimization models with decision variables in multiple time scales are widely used across various fields such as integrated planning and scheduling. To address scalability challenges in these models, we present the Parametric Autotuning Multi-Time Scale Optimization (PAMSO) algorithm. PAMSO tunes parameters in a low-fidelity model to assist in solving a higher-fidelity multi-time scale optimization model. These parameters represent the mismatch between the two models. PAMSO defines a black-box function with tunable parameters as inputs and multi-scale cost as output, optimized using Derivative-Free Optimization methods. This scalable algorithm allows optimal parameters from one problem to be transferred to similar problems. Case studies demonstrate its effectiveness on an MINLP model for integrated design and scheduling in a resource task network with around 67,000 variables and an MILP model for integrated planning and scheduling of electrified chemical plants and renewable resources with around 26 million variables.
\end{abstract}



\end{frontmatter}



\section{Introduction}

In large systems consisting of different interrelated processes, decisions are taken at different time scales including high-level planning decisions as well as low-level scheduling and control decisions. In the standard decision hierarchy adopted by the process systems engineering community,  decisions from high-level decision problems are passed on to low-level operating decision problems to obtain a complete execution \Citep{Shin2019Multi-timescaleProgramming}. High-level decision problems include planning problems (e.g., supply chain design, capacity planning, production targets, and task selection) and provide a set of decisions on a longer time scale, i.e., on a yearly, monthly, or weekly basis. Low-level operating decision problems include scheduling and control decisions and provide a set of decisions on a shorter time scale, i.e., on an hourly scale or in seconds. These decision layers are interconnected and, a one-way top-down communication in the hierarchical structure, i.e., solving the high-level problems without input from the low-level problems and then fixing the high-level decisions in the low-level problems, can result in infeasible or suboptimal solutions. Therefore, it is necessary to integrate the different decision layers. This can be achieved by solving multi-time scale decision problems which involve formulating a multi-time scale optimization model that integrates decision processes taking place at different time scales \citep{Shin2019Multi-timescaleProgramming,Maravelias2009IntegrationOpportunities,Allen2023ASystems,Subramanian2013IntegrationManagement}.  Multi-time scale optimization models are of significant importance in various domains, such as plant production and scheduling \citep{Biondi2017OptimizationApproach} and electricity market \citep{Dowling2017AParticipation} and environmental management\citep{ Tabrizi2018IntegratedImpacts}. In particular, an area where multi-time scale optimization plays a crucial role is in the decarbonization and electrification of chemical process systems \citep{Ramanujam2023DistributedMicrogrid,Kim2024Multi-periodUncertainty}. 

An important issue to deal with when using multi-time scale optimization models is the fact that they can have a very large number (in the order of tens of millions or more) of variables and constraints. While full-space methods \Citep{Kopanos2009Multi-SiteIndustries} have been used, it becomes difficult to solve these problems as they become larger (in the order of millions of variables). In order to solve these problems,  different methods such as decomposition methods \Citep{Erdirik-Dogan2008SimultaneousLines,Terrazas-Moreno2011APlants, Munawar2005IntegrationManufacturing,Shah2012IntegratedIndustry, Mora-Mariano2020AControl, Munoz2015SupplyEnvironment,Gharaei2019AScheduling,Peng2019AUncertainty,Barzanji2020DecompositionProblem,Chu2013IntegratedApproach}, metaheuristic methods \Citep{Lee2019SustainableRepresentation,Rauf2020IntegratedAlgorithm,Leite2023SolvingAlgorithms}, data-driven approaches \Citep{Dias2020IntegrationModels,Badejo2022IntegratingAnalysis,Beykal2022Data-drivenUncertainty,Wen2017AnSystem,Shin2019Multi-timescaleProgramming,Nikolopoulou2012HybridManagement,Ye2015AUncertainty,Yang2023IntegratedResources, Chu2014IntegratedModeling}, as well as mathheuristics \Citep{ Silva2023AProblem, Reinert2023ThisDecomposition, Ramanujam2023DistributedMicrogrid} have been applied. 

While decomposition methods can solve moderately large problems that have around millions of variables \citep{Li2022Mixed-integerSystems}, they face issues when applied to extremely large problems with more than tens of millions of variables, particularly when some of these variables are integers. Decomposition methods involve algorithms relying on iteratively solving a sequence of smaller problems until some convergence criteria are met. The rate of achieving convergence for extremely large problems could be a very slow process. Therefore, these algorithms are not scalable.  Furthermore, each of these decomposition algorithms is restricted to problems of a particular structure.  Scalability issues also occur in existing metaheuristic methods for problems with more than a thousand variables \citep{Hussain2019MetaheuristicSurvey} as well as matheuristic methods. Additionally, most existing data-driven methods are not scalable nor interpretable. Another important observation about these existing solution strategies is that these methods are not transferable. That is, solving a smaller version of the multi-time scale models using these methods does not help us solve a larger version of these models. Transfer learning is a technique in machine learning (ML) in which knowledge learned from one task is reused to boost performance on a related task \citep{Weiss2016ALearning}. Taking inspiration from transfer learning, a similar method for optimization problems can help solve problems faster.

In this paper, we propose the Parametric Autotuning Multi-time Scale Optimization algorithm (PAMSO) to deal with these issues.  Our algorithm is inspired by the parametric cost function approximations (CFAs) \Citep{PerkinsIII2017StochasticApproximations} from the reinforcement learning literature and the tuning of controllers using derivative-free optimizers  \Citep{Paulson2023ARepresentations, Coutinho2023BayesianControllers} from the control literature. Parametric CFA is a method that has been used to solve multistage sequential decision-making problems under uncertainty \Citep{Powell2022TheProgramming}.

The general idea behind the PAMSO algorithm is that we split the multi-time scale optimization model into a high-level decision layer and a low-level decision layer.  The high-level decision layer makes the high-level decisions, and the low-level decision layer takes into account the high-level decision variables and outputs the low-level decisions and the true objective we want to optimize. Without any further modifications, the mentioned idea refers to the one-way top-down communication and can lead to infeasibilities and bad solutions due to a mismatch between the high-level and low-level decision layers. To bridge this mismatch between the high-level and low-level decision layers, we add additional parameters to the high-level decision layer. We derive a good feasible solution to the problem by tuning these additional parameters. These tunable parameters reflect the uncertainties and detailed physics in the low-level model. Incorporating the parameters into the high-level decision layer helps improve the high-level decisions to obtain a better overall feasible solution.

For example, to solve the optimization model integrating the planning and scheduling of electrified chemical plants and renewable resources \citep{Ramanujam2023DistributedMicrogrid},  we can develop a simplified high-level model to output the planning decisions and fix the planning decisions in the integrated full-space model to obtain the scheduling decisions. The high-level model does not take into account all the detailed variations in the output of renewable resources. Without these detailed variations, some of the renewable resources may not be profitable in the high-level model but may be profitable in the integrated model. To remedy this, we can add a discount factor to the cost of the renewable resources in the high-level model, which serves as a tunable parameter, and tune the discount factor to make the renewable resources profitable in the high-level model.

We treat the entire hierarchy of the parameterized high-level decision layer and the low-level decision layer as a black box. The input of the black box is the set of tunable parameters and the output of the black box is the corresponding objective value of the multi-time scale model. We tune the parameters of the black box to obtain a good feasible solution.

We can incorporate transfer learning of the tunable parameters into the algorithm by choosing parameters that do not depend on the size of the problem. The algorithm can be implemented on a problem of small size to obtain parameters conducive to a good feasible solution. These parameters can then be used to obtain a good feasible solution for a similar problem with a much larger size in less time. More details of the algorithm are given in subsequent sections. The major contributions of this work are listed below:
\begin{itemize}
    \item We propose  Parametric Autotuning Multi-time Scale Optimization algorithm (PAMSO) as a scalable algorithm to solve multi-time scale optimization models.
    \item PAMSO can incorporate transfer learning into multi-time scale optimization models to help solve extremely large problems more quickly.
    \item PAMSO is applied to solve a set of case studies on integrated design and scheduling in a resource task network and integrated planning and scheduling of electrified chemical plants and renewable resources.
    \item PAMSO is available as an open-source software package with both implementation of the algorithm and case studies through our GitHub repository PAMSO.jl with URL: \url{https://github.com/li-group/PAMSO.jl}.
\end{itemize}

The rest of this paper is organized as follows. In Section \ref{sec:literaturereview}, a literature review of related works is provided. In Section \ref{sec:algorithm}, we describe the algorithm and its details. In Section \ref{sec:examples}, we illustrate the effectiveness of the proposed algorithm on a series of case studies. In Section \ref{sec:scalability}, we discuss the scalability of the proposed algorithm. The conclusions are drawn in Section \ref{sec:conclusion}.

\section{Literature Review} \label{sec:literaturereview}
\subsection{Algorithms for solving multi-time scale models}
 A variety of algorithms have been developed to solve multi-time scale optimization models. We classify these algorithms as follows:
\begin{enumerate}
\item Full space methods: In full space methods, a single optimization model is proposed for integrating the decisions layers, and that model is solved with exact optimization techniques. For example, \cite{Kopanos2009Multi-SiteIndustries} proposed an MILP to solve scheduling/batching and production planning for multi-site batch process industries and solved the MILP using the CPLEX solver. The issue with this method is that as the problem sizes increase, the problem becomes intractable. Using this method to solve problems with millions of variables and more is very time-consuming.
\item Decomposition methods: Decomposition methods can be used on problems moderately large (up to tens of millions of variables \Citep{Li2022Mixed-integerSystems}) and have a specific structure. The decomposition methods used for solving multi-time scale models include bi-level decomposition \Citep{Erdirik-Dogan2008SimultaneousLines,Terrazas-Moreno2011APlants}, multi-level decomposition \Citep{Munawar2005IntegrationManufacturing}, Lagrangian relaxation and/or decomposition \Citep{Shah2012IntegratedIndustry, Mora-Mariano2020AControl, Munoz2015SupplyEnvironment}, branch and price algorithm \Citep{Gharaei2019AScheduling}, progressive hedging \Citep{Peng2019AUncertainty}, and Benders decomposition method \Citep{Barzanji2020DecompositionProblem,Chu2013IntegratedApproach}.
\begin{enumerate}
    \item In bi-level decomposition, the original model is decomposed into a high-level and low-level model. For a multi-time scale problem, the high-level model can be a planning problem and the low-level problem can be a scheduling problem. The high-level model is usually a relaxation and provides bounds for the overall objective function and the low-level model gives a feasible solution. The high-level and low-level models are solved iteratively with problem-specific logic and integer cuts added to the high-level model based on the solution from the low-level problem. For example, \cite{Erdirik-Dogan2008SimultaneousLines} presented a multi-period MILP for the simultaneous planning and scheduling of single-stage multi-product continuous plants with parallel units and proposed a bi-level decomposition algorithm to solve the model. The multi-level decomposition algorithm involves extending the bi-level decomposition algorithm to have models at more than 2 levels. For example,  \cite{Munawar2005IntegrationManufacturing}  presented a multi-level decomposition-based framework for the integration of planning and scheduling in a multi-site, multi-product plant, with applications to paper manufacturing. 
\item Lagrangian decomposition (and relaxation) involves decomposing the original problem into smaller subproblems by dualizing the equalizing/complicating constraints. A Lagrangian master problem is solved to update the Lagrangian multipliers. Different extensions of Lagrangian decomposition/relaxation methods exist based on their implementation. For example, \cite{Shah2012IntegratedIndustry} used the augmented Lagrangian decomposition method to solve the integrated planning and scheduling problem for the multisite, multiproduct batch plants. Progressive hedging is a decomposition method for solving large-scale stochastic linear programs. The algorithm's core idea is to relax the non-anticipativity constraints and solve the scenario subproblems independently, and the algorithm penalizes deviations from the average values of decision variables. For example, \cite{Peng2019AUncertainty} proposed a novel Progressive hedging-based algorithm to address integrated planning and scheduling problems under demand uncertainty in a general mathematical formulation.
\item Branch and price is a hybrid algorithm of branch and bound and column generation algorithms.  Column generation algorithms involve using a subset of variables to solve the problem forming the restricted master problem. Variables are added iteratively to the subset by solving subproblems called pricing subproblems to enhance the solution \Citep{Desaulniers2006ColumnGeneration}. For example, \cite{Gharaei2019AScheduling}, used the branch and price algorithm to solve an integrated production scheduling and distribution problem with routing decisions in a multi-site supply chain.  
\item Benders decomposition is another decomposition algorithm used in solving multi-time scale optimization problems, where the problem is split into a high-level master problem and a set of low-level subproblems. The master problem and subproblems are solved iteratively,  with Benders cuts passed as feedback to the master from solving the subproblems. Classical Benders decomposition requires that the subproblems be linear \Citep{Chu2013IntegratedApproach}. There are different extensions of Benders decomposition like logic-based Benders decomposition \Citep{Barzanji2020DecompositionProblem} where the subproblems do not have to be linear. For example, \cite{Chu2013IntegratedApproach}  adapted the generalized Benders decomposition to solve a mixed-integer dynamic optimization integrating scheduling and dynamic optimization for batch chemical processes .\cite{Barzanji2020DecompositionProblem} implemented a logic-based Benders decomposition to solve an integrated process planning and scheduling problem.

\end{enumerate}
The performance of these decomposition techniques depends on the structure and size of the problem. For very large problems with say tens of millions of variables, using decomposition algorithms leads to a large number of iterations to converge which could require a large amount of time. Furthermore, these decomposition algorithms are restricted to problems of a particular structure.

\item Metaheuristic Algorithms:  A metaheuristic is a
high-level procedure or heuristic designed to find, generate, or select a heuristic (partial search algorithm) that may provide a sufficiently good solution to an optimization problem, especially with
incomplete or imperfect information. Examples include genetic algorithms \Citep{Lee2019SustainableRepresentation}, simulated annealing \Citep{Li2007AScheduling}, ant colony optimization \Citep{Liu2016ApplicationScheduling}, and variable neighborhood search \Citep{Leite2023SolvingAlgorithms}. For example, \cite{Lee2019SustainableRepresentation} used the genetic algorithm to find the pseudo-optimum of integrated process planning and scheduling problems. \cite{Leite2023SolvingAlgorithms} adapted the variable neighborhood search algorithm to address the integrated planning and scheduling problem on parallel and identical machines.  Metaheuristic methods are not scalable to problems to problems with more than a thousand variables \Citep{Hussain2019MetaheuristicSurvey}.

\item Data-driven methods: Surrogate models and feasibility analysis based on data have been used to describe the feasible regions of the scheduling models which can be then used to solve the integrated model \Citep{Maravelias2009IntegrationOpportunities, Dias2020IntegrationModels,Badejo2022IntegratingAnalysis, Ye2015AUncertainty,Yang2023IntegratedResources, Chu2014IntegratedModeling,Beykal2022Data-drivenUncertainty}. For example, \cite{Dias2020IntegrationModels} proposed a framework for the integration of planning, scheduling, and control using data-driven methodologies, which consisted of addressing the integrated problem as a grey-box optimization problem, and using data-driven feasibility analysis and surrogate models to approximate the unknown black box constraints.  \cite{Badejo2022IntegratingAnalysis} proposed a methodology for integrating scheduling operations into the supply chain network, motivated by the available enterprise data, and feasibility analysis. Data-driven methods based on other machine learning techniques \Citep{Wen2017AnSystem} have also been used.

Reinforcement learning-based approaches \Citep{Shin2019Multi-timescaleProgramming,Ochoa2022Multi-agentMarkets} have also been applied to integrated planning and scheduling optimization problems. Reinforcement learning algorithms, in particular, learn from interaction with the environment and can adapt to changing conditions over time. For example, \cite{Shin2019Multi-timescaleProgramming}, developed a multi-time scale decision-making model that combines Markov decision process (MDP) and mathematical programming (MP) in a complementary way and introduced a computationally tractable solution algorithm based on reinforcement learning (RL) to solve the MP-embedded MDP problem.  \cite{Ochoa2022Multi-agentMarkets}  proposed a novel multi-agent deep reinforcement learning framework for efficient multi-time scale bidding for hybrid power plants.

 Most of these methods are not scalable and not physically interpretable. 

\item  Matheuristics: Matheuristics \Citep{Ball2011HeuristicsProgramming} are algorithms that combine mathematical programming and heuristics. For example, \cite{Silva2023AProblem}, proposed a relax-and-fix heuristic to solve an integrated multiproduct, multiperiod, and multistage
capacitated lot sizing with a hybrid flow shop problem. \cite{Reinert2023ThisDecomposition} proposed a spatial aggregation and decomposition algorithm to solve large multi-time scale energy problems.  \cite{Ramanujam2023DistributedMicrogrid} proposed an aggregation-disaggregation based algorithm to solve an integrated planning and scheduling problem for electrified chemical plants. Matheuristics are efficient for problems with small to moderate size. For example, the matheuristic proposed by \cite{Ramanujam2023DistributedMicrogrid}, can solve a problem with around 1 million variables and 10 million constraints. However, matheuristics are not scalable to larger problems and are problem specific.
\end{enumerate}

\subsection{Parametric Cost Function}

\begin{figure}[!ht]
\centering
\begin{subfigure}{\textwidth}
\centering
\includegraphics[width=15cm]{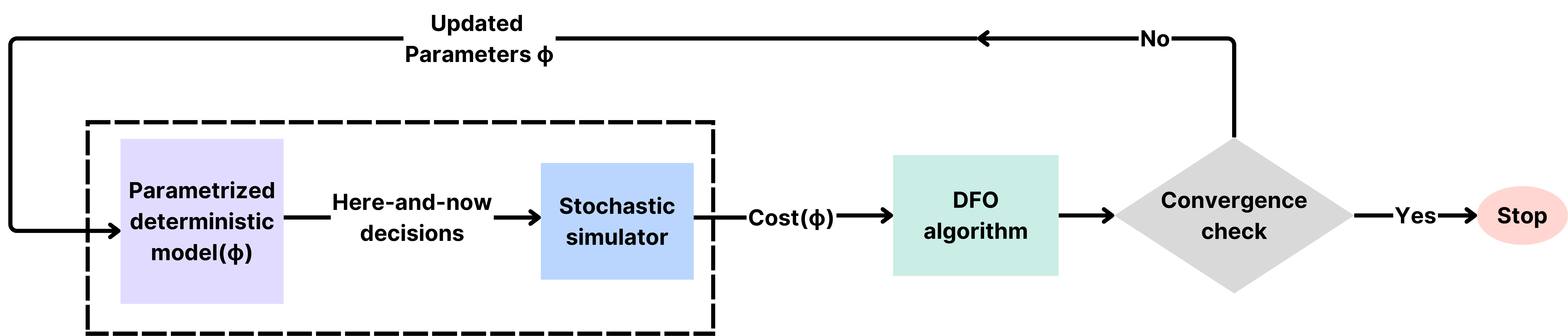}

\caption{Parametric CFA}
\label{fig:CFA}
\end{subfigure}

\bigskip
\begin{subfigure}{\textwidth}
\centering
\includegraphics[width=15cm]{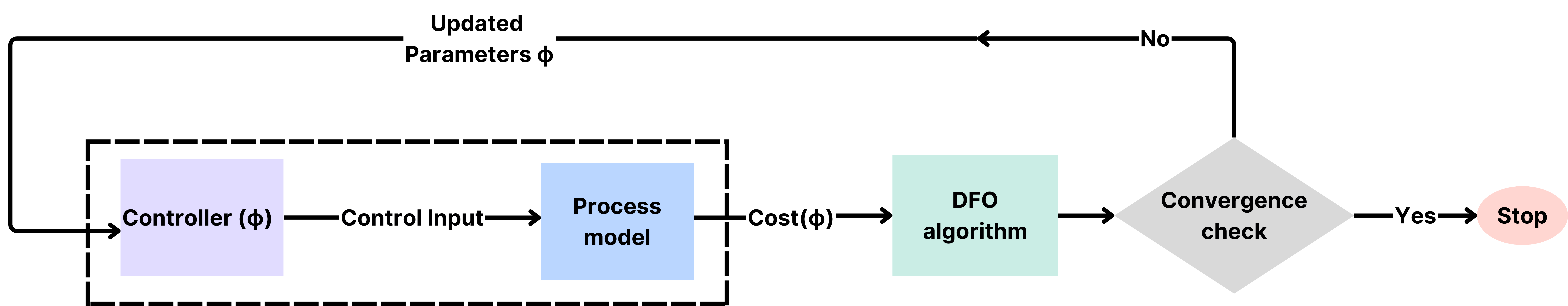}

\caption{Tuning of the controller through derivative-free optimization}
\label{fig:Controller}

\end{subfigure}
\begin{subfigure}{\textwidth}
\centering
\includegraphics[width=15cm]{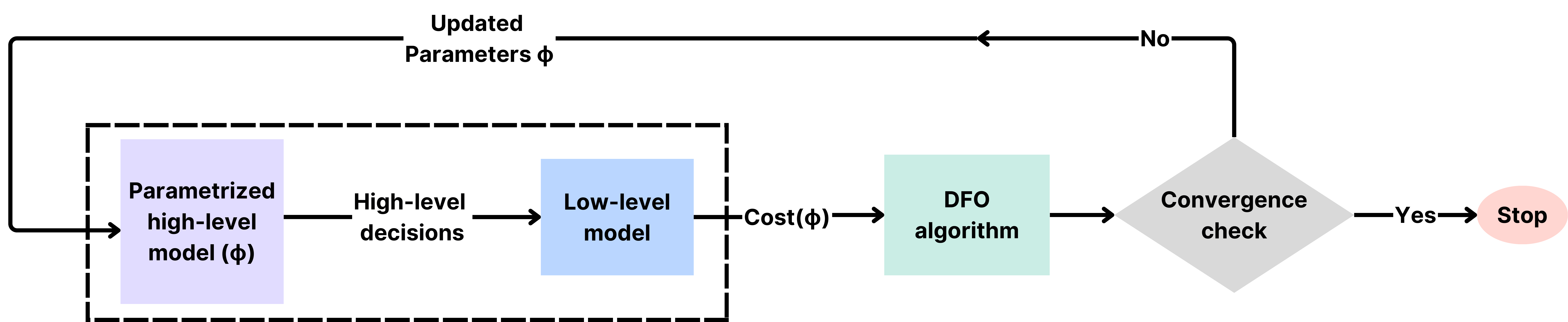}

\caption{PAMSO}
\label{fig:PAMSO}
\end{subfigure}

\caption{Comparison of parametric CFA, tuning of controller and PAMSO}

\end{figure}

We draw inspiration for the algorithm from the parametric cost function approximations (CFAs). CFAs are one of the strategies used for solving sequential decision-making problems under uncertainty. A common framework used by the PSE community to model sequential decision-making problems under uncertainty is the multistage stochastic programming, where the probability distribution of the uncertain parameters is approximated \Citep{Li2021AUncertainty}. The uncertainties are generally characterized by discrete realizations of the uncertain parameters as an approximation to the real probability distribution with each realization defined as a scenario. An issue with this framework is that the corresponding scenario tree increases exponentially with the number of stages making the corresponding optimization model extremely large, which causes scalablity issues even with the decomposition algorithms. Furthermore, stochastic programming models are approximations of real lookahead models and hence solving them may not provide the optimal solutions to the real problem.  By contrast, many industries solve deterministic optimization models 
to make sequential decisions. Deterministic optimization models are very fast to solve and easier to understand.  The flip side of using such deterministic optimization models is their inability to reflect uncertainty. To handle these uncertainties, there have been heuristics that involve tunable parameters. \cite{PerkinsIII2017StochasticApproximations} formally characterized these parametrized deterministic models as parametric cost function approximations (CFAs).  The idea behind parameterized CFA is to add some parameters into a low-fidelity deterministic model, such as the cost coefficients, and tune them using a high-fidelity ``stochastic simulator'' through derivative-free or derivative-based
optimization as shown in Figure \ref{fig:CFA}. CFAs are widely used for solving large single-scale problems such as planning a supply chain under uncertainty or airline scheduling  \Citep{Powell2019AOptimization}.

A simulation-based tuning approach is also utilized in tuning a controller with a derivative-free optimizer, as depicted in Figure 
\ref{fig:Controller} \Citep{Paulson2023ARepresentations, Coutinho2023BayesianControllers}. The control input is obtained from the parameterized low-fidelity controller, which is then passed on to the high-fidelity process model which acts as a simulator. The parameters of the controller are tuned based on the cost from the process model. The common theme in parameterized CFAs and controller tuning is to improve the performance of a parameterized low-fidelity model (parameterized deterministic model and controller model) by tuning the parameters using the feedback from the high-fidelity model (stochastic simulator and process model). The high-fidelity model acts like a simulator taking input from the parameterized low-fidelity model and outputting the performance of the entire system. 

Adapting this philosophy, we split the multi-time scale optimization model into a hierarchy of a high-level and low-level model. The high-level model is a low-fidelity model that outputs high-level decisions like planning decisions. The high-level model is analogous to the parametric deterministic model in CFAs and the controller model while tuning a controller. The low-level model acts as a high-fidelity simulator taking the high-level decisions as the input and outputting the low-level decisions and the true objective of the multi-time scale optimization model. The low-level model is analogous to the stochastic simulator in CFAs and the process model in the tuning of controllers. We parameterize the high-level model, with the parameters reflecting some of the mismatch between the high-level and low-level model. The parameterized high-level model and the low-level model are treated as a single system, and the parameters are tuned using derivative-free optimizers to optimize the net objective function of the multi-time scale model, as shown in Figure \ref{fig:PAMSO}.

 \section{Proposed Algorithm}\label{sec:algorithm}
In this section, we discuss the details of the algorithm. We sequentially split the decision variables in multiple time scales into two sets of variables: high-level decision variables denoted by $\bm{x}$ which pertain to the long time scales and low-level decision variables $\bm{y}_s$ encompassing the other shorter time scales defined for the scenarios/subperiods $s$.  
We write the multi-time scale optimization model as the following:
\begin{subequations}
    \label{FP}
    \begin{align}
       \min_{\bm{x},\bm{y}_s
       } \quad f(\bm{x}) + \sum_{s \in S}w_s q(\bm{x}, \bm{y}_s,\bm{\theta}_s) & \\
        \text{s.t.} \quad g(\bm{x}) \leq 0 & \\
        h(\bm{x}, \bm{y}_s, \bm{\theta}_s) \leq 0 \quad \forall s \in S & 
    \end{align}
\end{subequations}
Here, $\bm{\theta}_s$ represents the parameters associated with scenario/subperiod $s$ and \(w_s\) denotes the probability/weight of scenario/subperiod $s$. \(f\) represents the cost function associated with the high-level decisions. \(q\) denotes the cost function associated with the low-level decisions. \(g\) represents the constraints on the high-level decisions and \(h\) denotes the set of detailed constraints on the low-level decisions as well as the constraints connecting the high-level decisions with the low-level decisions.

To solve \eqref{FP}, we realize that a very fast way to get a solution is to split \eqref{FP} into a hierarchy of a high-level model \eqref{UP} and a low-level model \eqref{LP} and implement the one-way communication method on the hierarchy. The high-level model provides the high-level decisions. Fixing the high-level decisions in the low-level model provides the low-level decisions and the corresponding objective for \eqref{FP}. The high-level model does not take into account the detailed variations involved in the low-level model and is usually small making it fast to solve. The low-level model is also very fast to solve with the high-level decisions fixed. Thus this hierarchy of a high-level and a low-level model is very fast to implement.

The high-level model is shown below:
\begin{subequations}
    \label{UP}
    \begin{align}
        \min_{\bm{x},\bm{z}} \quad f(\bm{x})+Q(\bm{x},\bm{z},\bm{\Tilde{\theta}}) & \\
        \text{s.t.} \quad g(\bm{x}) \leq 0 & \\
        H(\bm{x}, \bm{z}, \bm{\Tilde{\theta}}) \leq 0 & 
    \end{align}
\end{subequations}

In this formulation, $\bm{z}$ acts as a surrogate for the low-level decisions, and $\bm{\Tilde{\theta}}$ acts as a surrogate parameter for the original parameters $\bm{\theta}_s$. Furthermore, \(Q\) and \(H\) are surrogates for the function \(q\) and the set of constraints \(h\), respectively i.e., the equivalent functions and constraints in terms of \(\bm{z}\) and $\bm{\Tilde{\theta}}$.

The dimensions of \(\bm{z}\) and $\bm{\Tilde{\theta}}$ can be much smaller than the total dimensions of each \(\bm{y_{s}}\) and \(\bm{\theta_{s}}\) respectively. For example, consider solving a problem with high-level decisions taken for a month and low-level decisions taken on an hourly basis. We propose a high-level model that takes decisions for the month removing the hourly complexity. We can sum up the hourly variables \(\bm{y}_{s,\tau}\) for each hour \(s \in \set{S}\) at each day \(\tau \in \set{D}\) of the month and use the summed-up variables as the surrogate for \(\bm{y}_{s,\tau}\) i.e.,  \(\bm{z} = \sum_{s \in \set{S}}\sum_{\tau \in \set{D}}\bm{y}_{s,\tau}\). The surrogate parameter $\bm{\Tilde{\theta}}$ could be the sum of the original parameters $\bm{\theta}_{s,\tau}$ i.e., \(\bm{\Tilde{\theta}} = \sum_{s \in \set{S}}\sum_{\tau \in \set{D}}{\bm{\theta}_{s,\tau}}.\)  For the example explained above, the dimensions of \(\bm{z} = \sum_{s \in \set{S}}\sum_{\tau \in \set{D}}\bm{y}_{s,\tau}\) and  \(\bm{\Tilde{\theta}} = \sum_{s \in \set{S}}\sum_{\tau \in \set{D}}{\bm{\theta}_{s,\tau}}.\) have been reduced by a factor of \(|\set{S}|\times |\set{D}|\).

The low-level model is shown below:
\begin{subequations}
    \label{LP}
    \begin{align}
        \min_{\bm{y}_s} \quad f(\bm{x}^*) + \sum_{s \in S}w_s q(\bm{x}^*, \bm{y}_s,\bm{\theta}_s) & \\ 
        \text{s.t.} \quad h(\bm{x}^*, \bm{y}_s, \bm{\theta}_s) \leq 0 \quad \forall s \in S & 
    \end{align}
\end{subequations}

In this case, $\bm{x}^*$ represents the high-level decisions obtained from solving the high-level model.

The resulting low-level problem after fixing  $\bm{x}$ to $\bm{x}^*$ may be suboptimal or infeasible. To address this issue, we draw inspiration from the parametric cost function approximation  (CFA) approach introduced by \cite{PerkinsIII2017StochasticApproximations} which uses parameterized deterministic optimization models to solve multi-stage stochastic programming problems. In the hierarchy we are dealing with, there is a mismatch between the high-level and low-level models due to the detailed physics-based variations embedded in the low-level model. To reflect these variations in the high-level model without making the high-level problem too large, we add tunable parameters to the high-level model. The parameterized high-level model is shown below:

\begin{subequations}
    \label{PUP}
    \begin{align}
        \min_{\bm{x},\bm{z}} \quad \Tilde{f}(\bm{x},\bm{\rho})+\Tilde{Q}(\bm{x},\bm{z},\bm{\Tilde{\theta}},\bm{\rho}) & \\
        \text{s.t.} \quad \Tilde{g}(\bm{x},\bm{\rho}) \leq 0 & \\
        \Tilde{H}(\bm{x}, \bm{z}, \bm{\Tilde{\theta}},\bm{\rho}) \leq 0 & 
    \end{align}
\end{subequations}
Here \(\bm{\rho}\) are the set of tunable parameters we introduce. \(f,Q,g\) and \(H\) are modified into \(\Tilde{f},\Tilde{Q},\Tilde{g}\) and \(\Tilde{H}\),
respectively taking into account the parameters \(\bm{\rho}.\)

We thus use the same hierarchy by obtaining $\bm{x}^*$ from \eqref{PUP} and fixing it in \eqref{LP} to determine the solution. Depending on the complexity and physics of the problem, we can fix the value of a subset of high-level decision variables $\bm{x}^*$ in \eqref{LP}, letting the other variables be optimized along with the low-level decision variables. We treat this entire hierarchy of \eqref{PUP} and \eqref{LP} as a system and obtain the best possible solutions by tuning the parameters \(\bm{\rho}.\) Finding the derivative of the objective from the system with respect to the parameters is very difficult, especially with the presence of binary variables. Therefore, we treat the system as a black box and define it as the Multi-time scale black box function \hyperref[MBBF]{(MBBF)}.  We use a Derivative-Free Optimization (DFO) solver to optimize the MBBF. The parameters chosen are significant to the physics of the problem and help bridge the mismatch between the high-level model and the low-level model. Tuning the parameters by solving the black box model yields parameters conducive to a good feasible solution. 

\setlength{\interspacetitleruled}{-.4pt}%
\begin{algorithm*}[t]
\label{MBBF}
  \SetKwFunction{FMain}{MBBF}
  \SetKwProg{Fn}{Function}{:}{}
  \Fn{\FMain{$\bm{\rho}$}}{
     \tcc{Solve the parameterized high-level model with parameters \(\bm{\rho}\)}
     \begin{equation*}
          x^*,z^* = \argmin1_{\bm{x},\bm{z}} \Tilde{f}(\bm{x},\bm{\rho})+\Tilde{Q}(\bm{x},\bm{z},\bm{\Tilde{\theta}},\bm{\rho}) \quad 
        \text{s.t.} \quad \Tilde{g}(\bm{x},\bm{\rho}) \leq 0, 
        \Tilde{H}(\bm{x}, \bm{z}, \bm{\Tilde{\theta}},\bm{\rho}) \leq 0 
     \end{equation*}
     \BlankLine
     \tcc{Solve the low-level model using the high-level model decisions obtained from the previous step to obtain the true objective value}
     \begin{equation*}
        y^*_{s} =  \argmin1_{\bm{y}_s}  f(\bm{x}^*) + \sum_{s \in S}w_s q(\bm{x}^*, \bm{y_s},\bm{\theta}_s) \quad
        \text{s.t.} \quad h(\bm{x}^*, \bm{y}_s, \bm{\theta}_s) \leq 0 \quad \forall s \in S 
        \end{equation*}
        \BlankLine
    \tcc{Return the corresponding objective value of the multi-time scale model }
    \Return{\(f(\bm{x}^*) + \sum_{s \in S}w_s q(\bm{x}^*, \bm{y}^*_s,\bm{\theta}_s)\)}
  }

\end{algorithm*}

An important challenge is the feasibility of parameters in the multi-time scale black box function. There can be cases where the parameters yield infeasible high-level models or when the generated high-level decisions are infeasible in the low-level model. To handle these infeasibilities, the approach we take is to return a large positive value (eg., 1e10) from the multi-time scale simulator. Most DFO solvers balance the trade-off between exploration and exploitation. If
the initial set of parameters is infeasible, we increase the exploration in the DFO solver. 

While tuning the parameters, we can set a time limit and an MIP gap or any other related metric to the optimization models to ensure that each iteration is completed within a set period of time. After obtaining the best set of parameters, we can run the low-level model to optimality. Furthermore, based on the physics of the problem, we can decompose the low-level model into a set of optimization sub-models, such as by subperiod 
\(s\), whenever feasible. This decomposition helps solve the low-level model faster.

\textbf{Transferability of the tunable parameters:} To make the optimal tunable parameters transferable from one problem to similar problems, we can select parameters that remain independent of the problem size or independent of that part of the problem that changes.  We can then pre-train the parameters for a problem in one class and further fine-tune these trained parameters to solve a slightly different problem. For example, we can solve a problem small in size and obtain the optimal parameters. We can then use these optimal parameters as a starting point to solve a similar larger problem while constraining the search space to a smaller range centered around the starting point.

\section{Computational results} \label{sec:examples}
To test the proposed algorithm, we implement the algorithm on a series of case studies on integrated design and scheduling in a resource task network \citep{Castro2005SimultaneousFormulations, Barbosa-Povoa1999DesignAlgorithm} and integrated planning and scheduling of electrified chemical plants and renewable resources \Citep{Ramanujam2023DistributedMicrogrid}. The code for the examples can be accessed at 
\url{https://github.com/li-group/PAMSO.jl}. 

For the case studies we tackle, the proposed algorithm and the corresponding models are implemented in Julia/JuMP. All the MILPs and MINLPs are solved using Gurobi version 11.0 on a Linux cluster with 48 AMD EPYC 7643 2.3GHz CPUs and 1 TB RAM using 8 threads. The three DFO algorithms we consider applying to optimize the MBBF are the Mesh Adaptive direct search (MADS) (\cite{Audet2006MeshOptimization}), Bayesian optimization (BayesOpt) (\cite{Paulson2023ARepresentations,Martinez-Cantin2014BayesOpt:Bandits}), and particle swarm optimization (PSO) (\cite{Mogensen2018Optim:Julia}). The implementations of the DFO algorithms are accessed through NOMAD.jl\repeatedfootnote{1}{https://github.com/bbopt/NOMAD.jl}, BayesOpt.jl\repeatedfootnote{2}{https://github.com/jbrea/BayesOpt.jl}, and Optim.jl\repeatedfootnote{3}{https://github.com/JuliaNLSolvers/Optim.jl/} respectively.

\subsection{Integrated Design and Scheduling in a Resource Task Network} 

In this example, we implement the algorithm on the integrated design and scheduling for the production of chemicals in a resource task network (RTN). Resource task network (RTN) \citep{Pantelides1994UnifiedScheduling} is a scheduling model used by the PSE community  that treats equipment, materials, and utilities as resources and tasks as operations to transform these resources. The equipment used needs to be sized appropriately, considering both the demand for materials/chemicals and the scheduling dynamics of the process. To maximize the economic potential of the network, it is essential to integrate the design and scheduling decisions for the reaction network. To optimize the integration, we formulate a multi-time scale optimization model by extending the RTN model to an integrated design and scheduling model. Previous models have been proposed for the integration of design and scheduling in an RTN with some formulations using MILP approximations \citep{Castro2005SimultaneousFormulations, Barbosa-Povoa1999DesignAlgorithm}. The full-space model we use is shown in Appendix \ref{rtnformulation}.

\textbf{Problem Statement:} We have a network of chemical reactions where feed materials are converted to desired products through intermediates. The demand profile for the products consisting of the total demand for each product on each day is known. We also know the various parameters of the chemical reactions in the network including the interaction parameters between the tasks and resources in the network. The processes take place in vessels that need to be designed. The materials are stored in storage vessels, which also have to be designed. We assume that the reaction network starts with no inventory. 

Additionally, we also need to schedule the chemical reactions and determine the corresponding amount of production for each hour.  The scheduling or amount of production depends on the demand profile.  If the entire demand is not met, a penalty based on the unsatisfied demand is incurred. While the buying and selling of materials takes place on an hourly basis, the total demand for a day needs to be satisfied.

The goal is to maximize the overall profit which is the difference between the revenue and cost. The cost includes the amortized costs of the storage and reaction vessels, the cost of feed, the penalty paid, and the cost of implementing tasks. The revenue is the money obtained from selling the products. The prices of the materials are on a per-unit basis with the cost of transportation included.

\textbf{Strategy:} As mentioned before, we formulate a multi-time scale optimization model to integrate the design and scheduling decisions. The optimization model takes decisions in 3 time scales: single-time design decisions, daily decisions involving the unfulfilled demand, and hourly decisions involving the amount of materials bought and sold, the inventory of the materials as well as the scheduling and batch sizes of tasks. Sizing of vessels involves economies of scale, i.e., larger vessels often have a lower cost per unit volume compared to smaller vessels. We include the effects of economies of scale in our optimization model making it a Mixed Integer Non-Linear Programming (MINLP) model. 

The MINLP model can be very large with tens of thousands of variables and can be intractable with state-of-the-art algorithms. By implementing PAMSO we can obtain good feasible solutions within a reasonable amount of time. While implementing PAMSO on the integrated problem we use a high-level model that involves decisions in 2 time scales: single-time design decisions and daily decisions. The daily decisions involve the decisions on unfulfilled demand as well as the daily operating decisions. The daily operating decision variables are constructed by summing up the hourly operating decisions for each day. The high-level model is formulated using different methods of aggregation of the original demand profile which will be explored in detail in Sections \ref{1weekrtn} and \ref{4weekrtn}. We then pass the design decisions we obtain from the high-level model to the low-level model which gives the optimized scheduling decisions. The low-level model takes decisions in 2 time scales: daily decisions on unfulfilled demand, and hourly operating decisions involving the scheduling and batch sizes of tasks as well as the buying, selling, and inventory of materials. 
 
 In order to improve the performance of the high-level model to produce good design decisions, we modify the objective function appropriately by choosing tunable parameters to modify the weights to the different components of the objective function. This can have different effects on the design decisions. For example, increasing the cost of the penalty for unsatisfied demand and decreasing the cost of the vessels will lead to a design that tries to fulfill more of the demand at the cost of larger-sized vessels. On the other hand, increasing the cost of the vessels, and decreasing the penalty on unsatisfied demand or the prefactor on startup costs of tasks, will reduce the size of the vessels, but have more unfulfilled demand. We optimize the weights on different components of the objective to get an optimum design with the best total cost. To do this, we add tunable parameters to modify the penalty paid for unsatisfied demand, the cost of reaction and storage vessels as well as the cost of starting up tasks. We further classify the reaction vessels into three categories to capture more information about the network: vessels with at least one feed-based task (feed task vessels), vessels with no feed-based tasks but at least one product-based task (product task vessels), and vessels involving only intermediate tasks (intermediate task vessels). We assign tunable parameters to each of these types of vessels. A summary of the tunable parameters used is shown in Table \ref{table:rtnparameters}.The upper bound for the parameters was chosen empirically.

\begin{table}[H]
\centering
\caption{Parameters for implementing PAMSO on integrated design and scheduling model}
\label{table:rtnparameters}
\begin{tabularx}{\textwidth}{>{\centering\arraybackslash}X>{\centering\arraybackslash}X>{\centering\arraybackslash}X>
{\centering\arraybackslash}X}
\hline
Notation &  Parameter & Minimum value & Maximum value\\
\hline
$\rho_1$ & Prefactor to penalty for demand unsatisfied  & 0 & 30 \\
$\rho_2$ & Prefactor to cost of feed task vessels  & 0 & 30 \\
$\rho_3$ & Prefactor to cost of product task vessels  & 0 & 30\\
$\rho_4$ & Prefactor to cost of other intermediate task vessels  & 0 & 30  \\
$\rho_5$ & Prefactor to cost of storage units &  0 & 30 \\
$\rho_6$ & Prefactor to penalty for startup cost &  0 & 50 \\
\hline
\end{tabularx}
\end{table}

\begin{figure}[!ht]
\centering
 \begin{subfigure}[b]{0.4\textwidth}
 \includegraphics[width=8cm]{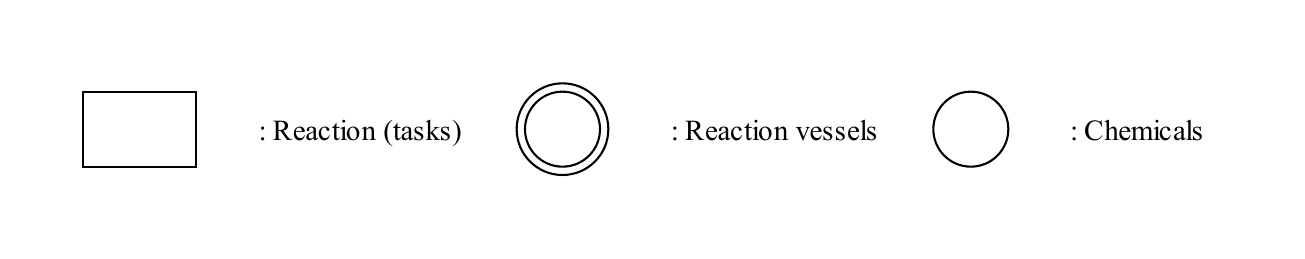}
 \end{subfigure}
 
 \begin{subfigure}[b]{0.8\textwidth}
\includegraphics[width=15cm]{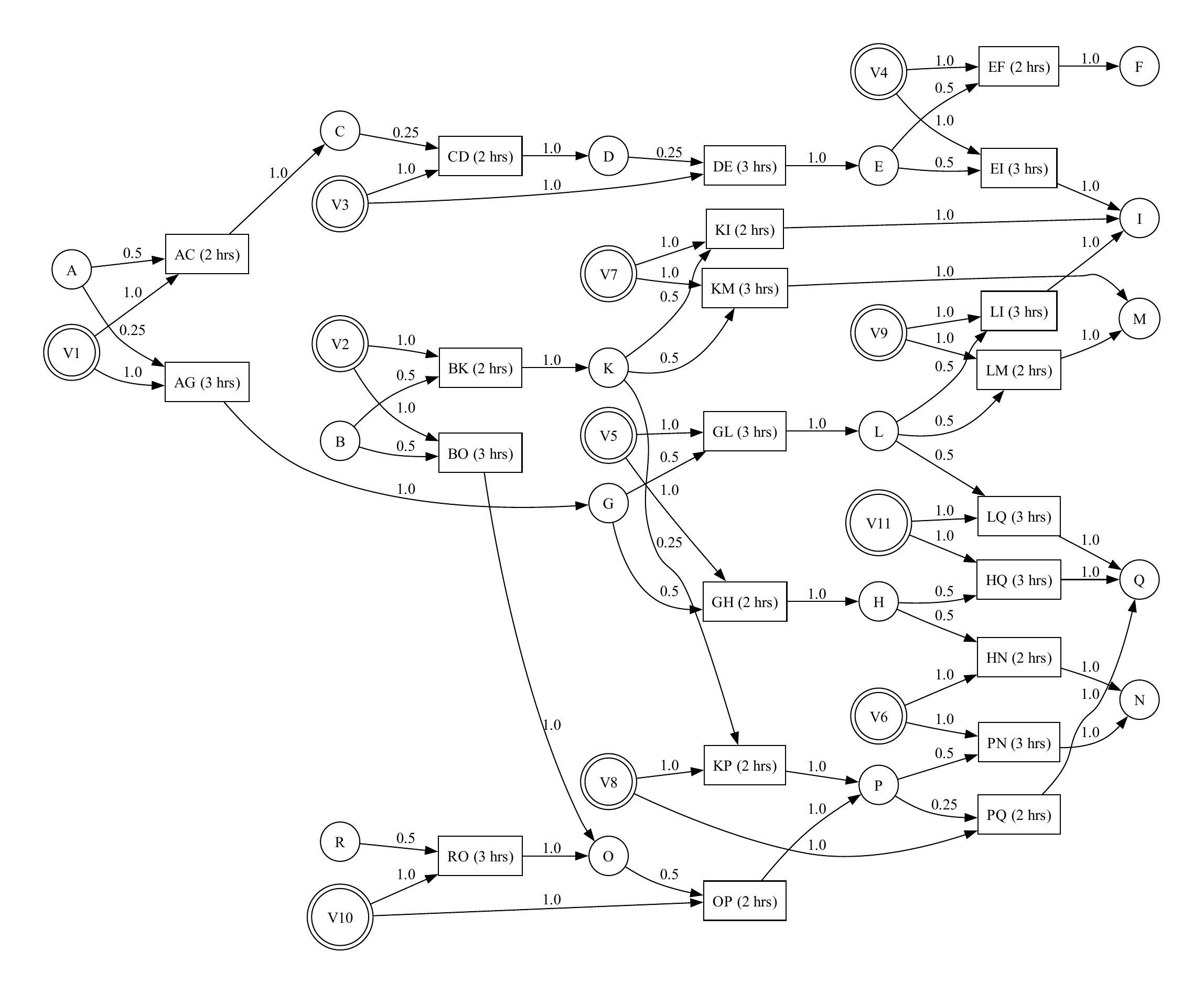}
\end{subfigure}
\vfill

\caption{Representation of RTN network}
\label{fig:repprocess}
\end{figure}

\textbf{Case study:} We consider an example of a network with 22 processes and 17 materials adapted from the example in \url{https://github.com/JavalVyas2000/rtn_scheduling/tree/main} with a few modifications. The 22 processes can occur in 11 vessels which are designed.  Additionally, each of the 17 materials has a storage vessel of its own which are also designed. The maximum size of a storage vessel is 100 units. The feed materials are bought and then converted into the products (F,I,M,N,Q) which are sold.  A graphical representation of the system is shown in Figure \ref{fig:repprocess}.

We first optimize the system using a demand profile for one representative week, considering daily demand variations. We then optimize the average profit obtained from the system based on a daily demand profile for 4 representative weeks.

\subsubsection{Problem with 1 representative week} \label{1weekrtn}
\begin{figure}[!ht]
\centering
\includegraphics[width=10cm]{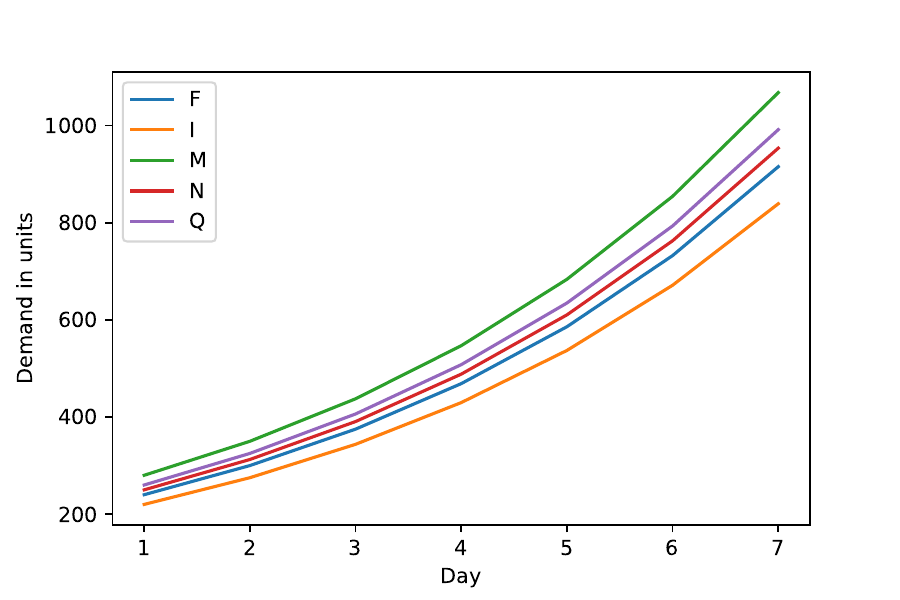}
\caption{Demand profile for representative week}
\label{fig:Dem_1week}
\end{figure}

In this part, we optimize the design and scheduling decisions for one representative week with a known demand profile. We implement PAMSO on the integrated problem. The demand profile for the representative week is used to formulate the high-level model. The low-level model is the integrated model with the design decisions fixed and takes only scheduling decisions and thus is an MILP model. The demand profile is shown in Figure \ref{fig:Dem_1week}. 

We apply the MADS, Bayesopt, and PSO algorithms to solve the associated MBBF for 150 function evaluations each. We start with \(\rho_1 = 1, \rho_2 = 0, \rho_3 = 0,\rho_4 = 0, \rho_5 = 0, \rho_6 = 0\) while applying the MADS and PSO algorithms. BayesOpt chooses its own initial parameters. The convergence of the algorithms for the case study is shown in Figure \ref{fig:iter_1week}. The figure shows that  PSO provides the best solution compared to the other two DFO solvers at the end of 150 function evaluations. 
The computing time and the optimal parameters are shown in Table \ref{table:rtn1weekres}. The best profit we obtain after applying the DFO solvers is \$8066.86.  From the set of optimal parameters, we notice that parameters modifying the penalty for the unfulfilled demand and startup cost are higher compared to the other parameters. This implies that for the system under consideration, it is very profitable to focus on minimizing the unsatisfied demand as compared to minimizing the size of vessels. This fact is further reiterated by the solution obtained from using the best parameters as the solution satisfies the entire demand of the customer. Another aspect to take into account is that the system is degenerate with respect to the parameters i.e., different sets of parameters can provide similar solutions. For example, while running the DFO solvers, we encountered a set of parameters, \(\rho_1 = 16.87,\rho_2 = 11.95,\rho_3 = 8.27,\rho_4 = 10.56, \rho_5 = 0.10, \rho_6 = 8.44\) which gave a profit of  \$6523.71.  Another set of parameters encountered, \(\rho_1 = 30, \rho_2 = 13.69, \rho_3 = 20.99, \rho_4 = 23.63, \rho_5 = 0.33, \rho_6 = 43.18\) gave a profit of \$6516.90. We can see that though the two sets of parameters are quite different, they yield similar profits. This shows the degenerate nature of the parameters due to their competing effects as well as the varying scales of the cost components. DFO solvers usually do not explore all possible directions equally. Hence, there is a chance of getting the best profit we obtained with a different set of parameters.

The optimal vessel size and storage size are shown in Tables \ref{table:vesresults_1week} and \ref{table:storresults_1week}, respectively. From the tables, we can see that some vessels have a size of 0 units, indicating that those vessels do not have to be installed for implementing the network. The corresponding reactions for the reaction vessels of size 0 units are not implemented to optimize the network. Additionally, some of the materials do not have a storage vessel installed as they are immediately consumed and/or produced by the network and hence do not require storage. 
 
 We also attempt to solve the full-space problem using Gurobi. The full space model has 13,223 continuous variables, 3,696 binary variables, and 23,583 constraints. The best profit we obtain after 24 hours by solving the full-space model using Gurobi is -\$8569.38 with an MIP gap of 850\%. Thus PAMSO outperforms Gurobi by giving a much better solution in a small amount of time.

\begin{figure}[!ht]
\centering
\includegraphics[width=10cm]{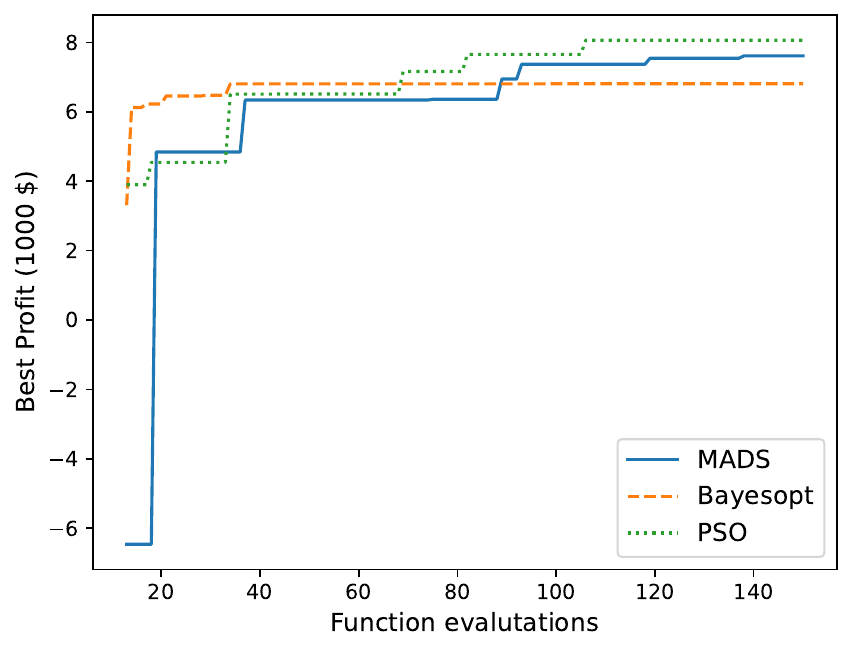}
\caption{Plot for convergence for 1 representative week}
\label{fig:iter_1week}
\end{figure}

\begin{table}[H]
\centering
\caption{Performance of DFO Algorithm  for 1 week}
\label{table:rtn1weekres}
\begin{tabular}{llll}
\hline
                 & MADS     & PSO      & Bayesopt \\
\hline
Time (h) & 3.45 & 4.61 & 3.61 \\
$\rho_1$ & 29.62 & 30 & 22.92 \\
$\rho_2$ & 29.18 & 17.84 & 11.95 \\
$\rho_3$ & 28.60 & 28.16 & 7.51 \\
$\rho_4$ & 16.60 & 27.60 & 10.55 \\
$\rho_5$ & 2.30 & 0.59 & 0.06 \\
$\rho_6$ & 43.90 & 50 & 8.30 \\
Best Profit (\$) & 7613.61 & 8066.86 & 6810.88 \\
\hline
\end{tabular}
\end{table}

\begin{table}[!ht]
\centering
\caption{Reaction vessel size for 1 representative week}
\label{table:vesresults_1week}
\begin{tabular}{ll}
\hline
Vessel & Size (units)\\
\hline
V1 & 2.63 \\
V2 & 72.58 \\
V3 & 48.35 \\
V4 & 67.96 \\
V5 & 1.50 \\
V6 & 95.37 \\
V7 & 184.83 \\
V8 & 112.35 \\
V9 & 0 \\
V10 & 0 \\
V11 & 0 \\
\hline
\end{tabular}
\end{table}

\begin{table}[H]
\centering
\caption{Storage size for 1 representative week}
\label{table:storresults_1week}
\begin{tabular}{ll}
\hline
Storage & Size (units)  \\
\hline 
A & 0 \\
B & 0 \\
C & 11.49 \\
D & 72.52 \\
E & 100 \\
F & 100 \\
G & 4.89 \\
H & 100 \\
I & 100 \\
K & 100 \\
L & 0 \\
M & 100 \\
N & 100 \\
O & 0 \\
P & 100 \\
Q & 100 \\
R & 0 \\
\hline
\end{tabular}
\end{table}
\subsubsection{Problem with 4 representative weeks}\label{4weekrtn}
In this part, we maximize the average profit obtained from designing and scheduling the network for 4 different demand profiles corresponding to 4 representative weeks. The demand profiles are shown in Figure \ref{fig:dem4rep}. We apply PAMSO to solve the problem.

\begin{figure}[htbp]
\captionsetup[subfigure]{justification=centering}
    \centering
    \begin{subfigure}[b]{0.48\textwidth}
        \centering
    \includegraphics[width=1\linewidth]{Images/Dem_week1.pdf}
\caption{Demand profile for representative week 1 }

    \end{subfigure}%
\hfill
    \begin{subfigure}[b]{0.48\textwidth}
        \centering
    \includegraphics[width=1\linewidth]{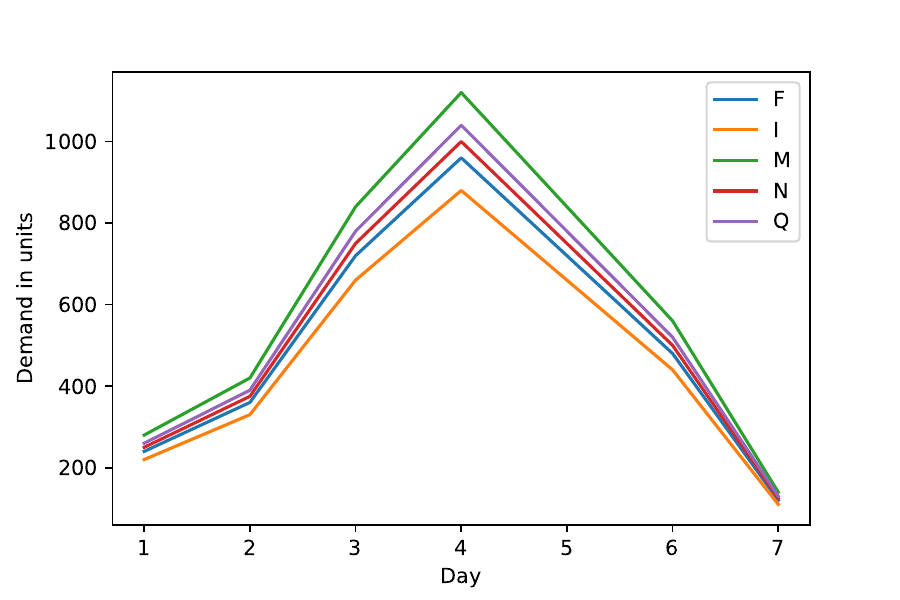}
\caption{Demand profile for representative week 2 }

    \end{subfigure}
\vspace{1em}
\captionsetup[subfigure]{justification=centering}
\ContinuedFloat 
\centering

    \begin{subfigure}[b]{0.48\textwidth}
        \centering
   \includegraphics[width=1\linewidth]{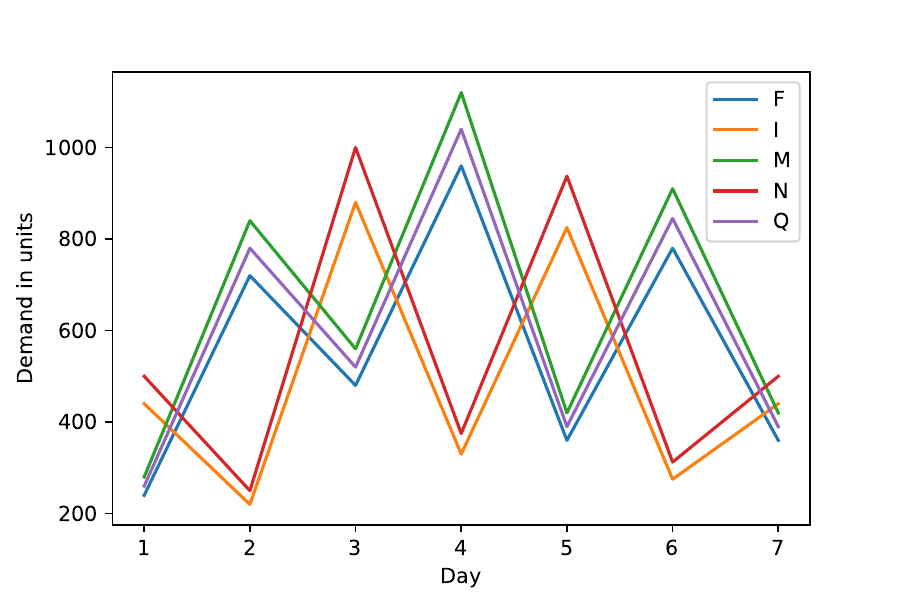}
\caption{Demand profile for representative week 3 }

    \end{subfigure}
\hfill
    \begin{subfigure}[b]{0.48\textwidth}
        \centering
\includegraphics[width=1\linewidth]{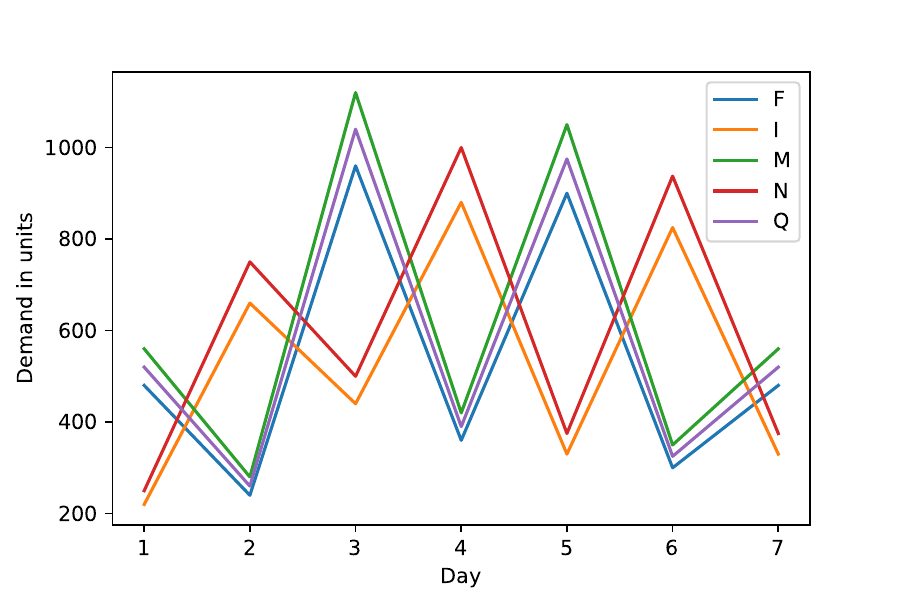}
\caption{Demand profile for representative week 4 }

    \end{subfigure}
   
    \caption{Demand profile for the 4 representative weeks}
 \label{fig:dem4rep}
\end{figure}

While implementing PAMSO, we apply the MADS, Bayesopt, and PSO algorithms to solve the associated MBBF for 150 function evaluations each. We start with \(\rho_1 = 1, \rho_2 = 0, \rho_3 = 0,\rho_4 = 0, \rho_5 = 0, \rho_6 = 0\) while applying the MADS and PSO algorithms. BayesOpt chooses its own initial parameters. 

\begin{figure}[H]
\centering
\includegraphics[width=10cm]{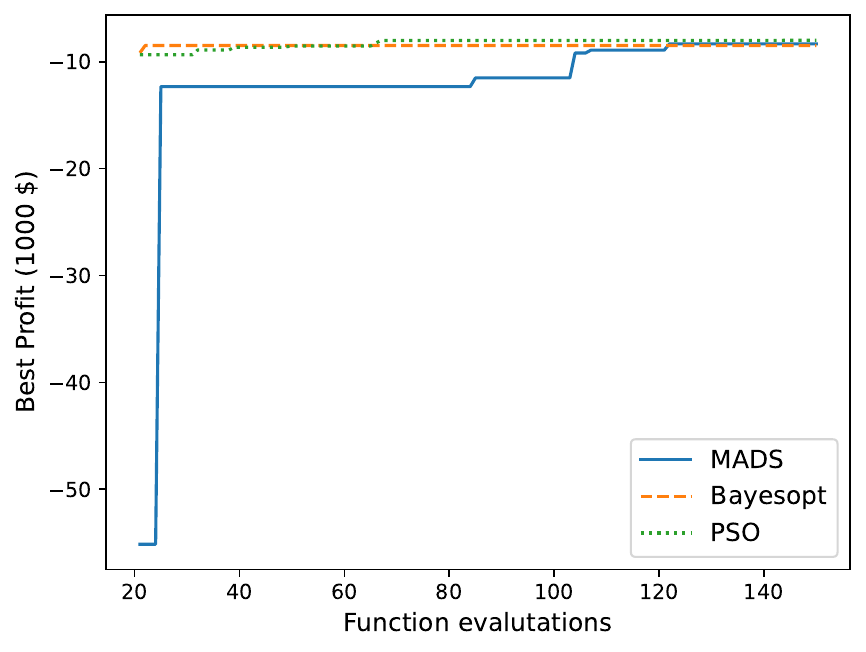}
\caption{Plot for convergence for approach 1}
\label{fig:iter_4week_agg}
\end{figure}

\begin{figure}[H]
\centering
\includegraphics[width=10cm]{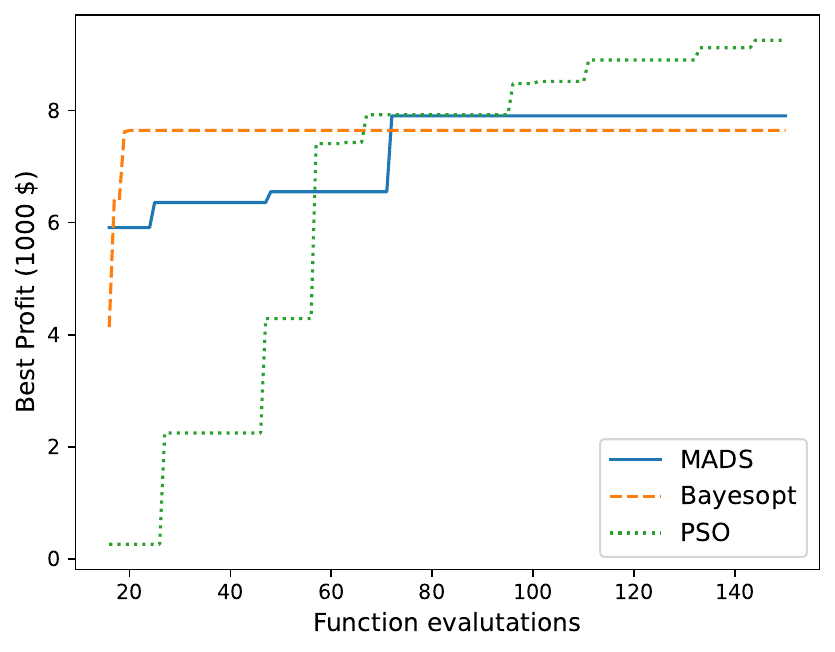}
\caption{Plot for convergence for approach 2}
\label{fig:iter_4week_dis}
\end{figure}

\begin{figure}[H]
\captionsetup[subfigure]{justification=centering}
    \centering
    \begin{subfigure}[b]{0.48\textwidth}
        \centering
    \includegraphics[width=1\linewidth]{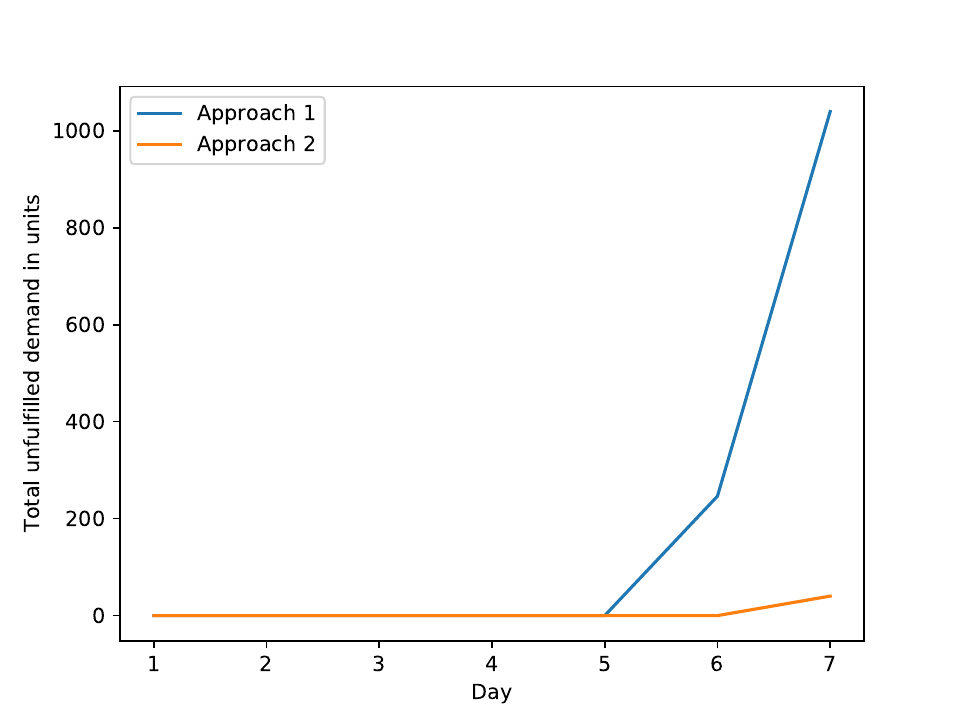}
\caption{Total unfulfilled demand for representative week 1 }

    \end{subfigure}%
\hfill
    \begin{subfigure}[b]{0.48\textwidth}
        \centering
    \includegraphics[width=1\linewidth]{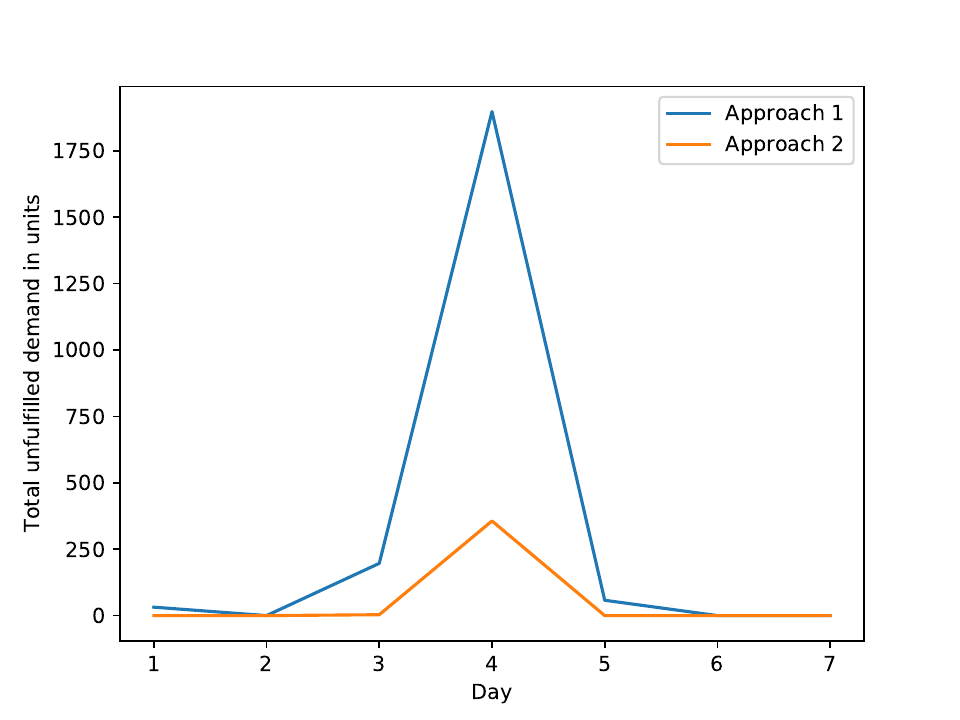}
\caption{Total unfulfilled demand for representative week 2 }

    \end{subfigure}
\vspace{1em}

    \begin{subfigure}[b]{0.48\textwidth}
        \centering
   \includegraphics[width=1\linewidth]{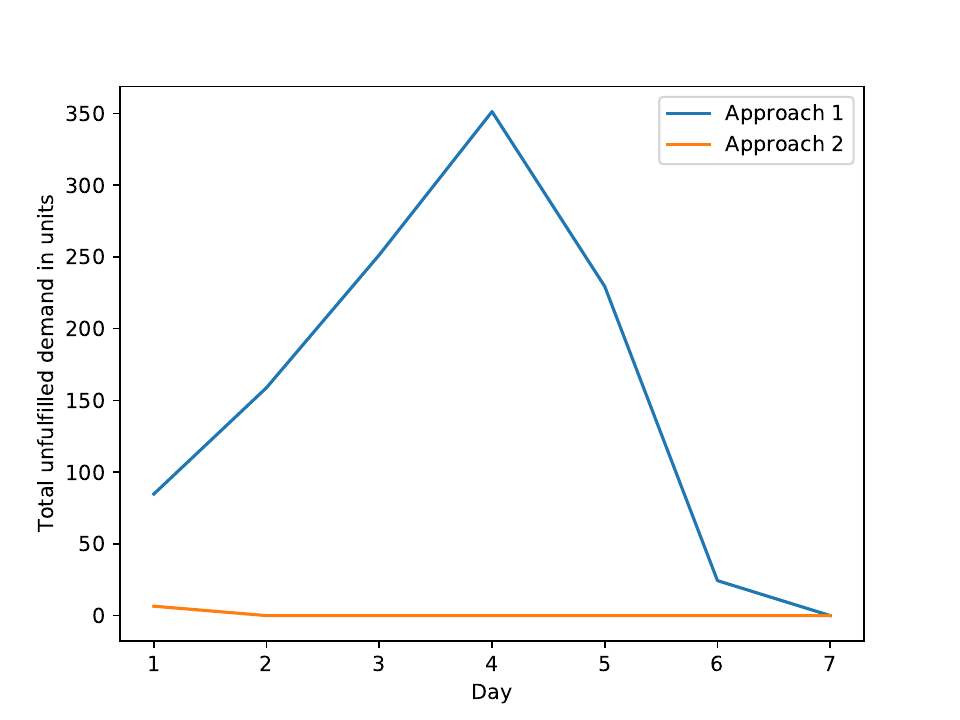}
\caption{Total unfulfilled demand for representative week 3 }

    \end{subfigure}
\hfill
    \begin{subfigure}[b]{0.48\textwidth}
        \centering
\includegraphics[width=1\linewidth]{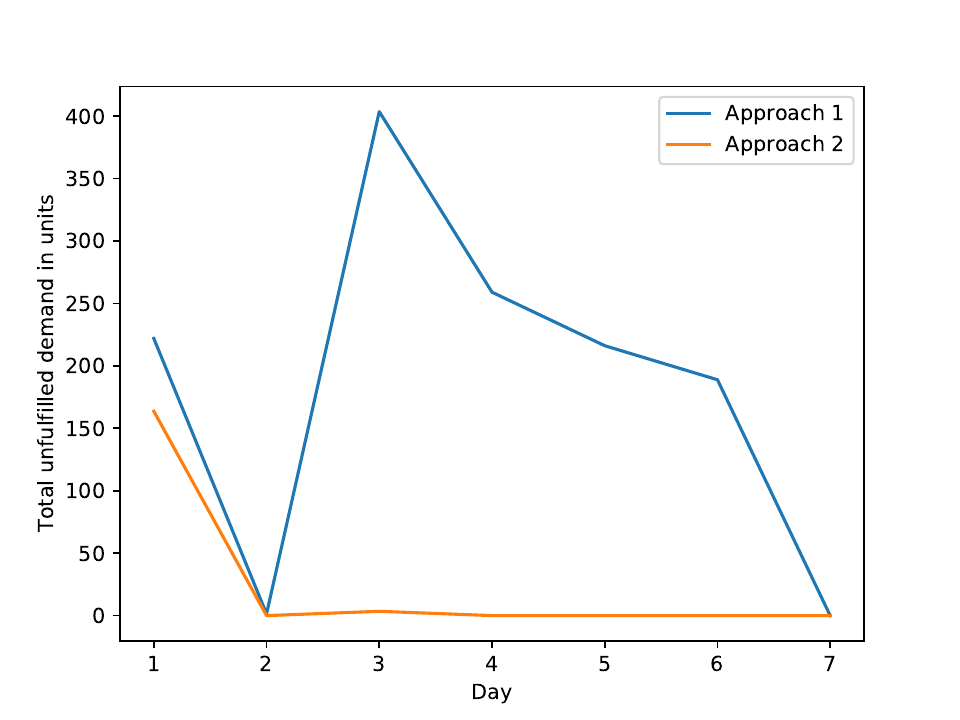}
\caption{Total unfulfilled demand for representative week 4 }

    \end{subfigure}
   
    \caption{Total unfulfilled demand for the 4 representative weeks}
 \label{fig:slack4rep}
\end{figure}

We can approach the problem in two ways: either by forming a high-level problem that involves the average demand profiles for the four weeks or by creating a more detailed high-level model that explicitly involves the demand profiles for each of the four representative weeks. Using the average of the demand profiles results in a smaller, less detailed high-level model, which yields a faster but potentially less accurate solution compared to the more detailed high-level model. The two approaches to s the problem are more formally defined as follows :

\begin{enumerate}
    \item \textbf{Approach 1}: We take the average daily demand for all the representative weeks and construct a demand profile over 7 days. More specifically, let \(D_{i,j,c}\) be the demand for chemical \(c\) on the \(j^{th}\) day of representative week $i$. The demand profile \(\Tilde{D}_{j,c} = \frac{\sum_{i=1}^4{D_{i,j,c}}}{4}\) is the demand profile we construct which represents the demand for chemical \(c\) on the \(j^{th}\) day. This demand profile over 7 days is used as the demand profile to formulate the high-level model. The low-level model which is the full-space model with design decisions fixed is decomposed by week. In other words, the low-level model consists of the set of integrated models for each representative week with the design decision fixed to those obtained from the high-level model. Solving this set of models will give the scheduling decisions for each representative week. The average of the profit obtained from solving this set of models is optimized while implementing PAMSO. The convergence of the algorithms for the case study is shown in Figure \ref{fig:iter_4week_agg}. From the figure, we can see that PSO gives the best solution as compared to the other 2 DFO solvers at the end of 150 function evaluations. 
The computing time and the optimal parameters are shown in Table \ref{table:rtn4weekres_agg}.  
 The best average profit we obtain after applying the DFO solvers is -\$7987.06. 
     \item \textbf{Approach 2}: We take the daily demand for all the representative weeks and construct a demand profile involving 28 days to solve the high-level model. More specifically, we connect the days of the representative weeks sequentially, to form a 28-day demand profile. The demand of the first 7 days is the demand from representative week 1, the demand from days 8 to 14, are the demand from representative week 2 and so on. This demand profile over 28 days is used as the demand profile to formulate the high-level model. The inventory at the end of each set of 7 days is fixed to 0.  The low-level model which is the full-space model with design decisions fixed is decomposed by week similar to that in approach 1. In other words, the low-level model consists of the set of integrated models for each representative week with the design decision fixed to those obtained from the high-level model. Solving this set of models will give the scheduling decisions for each representative week. The average of the profit obtained from solving this set of models is optimized while implementing PAMSO.  We use the tunable parameters defined in Table \ref{table:rtnparameters}. The convergence of the algorithms for the case study is shown in Figure \ref{fig:iter_4week_dis}. From the figure, we can see that PSO gives the best solution as compared to the other 2 DFO solvers at the end of 150 function evaluations.  The computing time and the optimal parameters are shown in Table \ref{table:rtn4weekres_dis}.  
 The best average profit we obtain after applying the DFO solvers is \$9254.19.  
\end{enumerate}
 The design decisions for the 2 approaches are shown in Table \ref{table:vesresults_4week} and Table \ref{table:storresults_4week}. We can see that the optimal vessel sizes in approach 2 are larger than the optimal vessel sizes in approach 1 ensuring that there is lesser unfulfilled demand in approach 2 and hence more profit. We also plot the total unfulfilled demand for each day in Figure \ref{fig:slack4rep} which reflects this trend. This is further emphasized by the best parameters as shown in Tables \ref{table:rtn4weekres_agg} and \ref{table:rtn4weekres_dis}. The best parameters in approach 1 are higher for the cost of feed task and intermediate task vessels than those in approach 2, ensuring that vessel sizes are smaller in approach 1. Furthermore, on a relative scale, the parameter for the penalty in approach 2 on the unfulfilled demand is higher than the other parameters ensuring that the unfulfilled demand is less in the solution. 

\begin{table}[!ht]
\centering
\caption{Reaction vessel size for 4 representative weeks}
\label{table:vesresults_4week}
\begin{tabular}{llll}
\hline
Vessel    & 4-week (approach 1) size (units) & 4-week (approach 2) size (units) \\
\hline
V1 & 1.87 & 2.21 \\
V2 & 60.48 & 70.72 \\
V3 & 43.01 & 52.73 \\
V4 & 56.28 & 71.67 \\
V5 & 0 & 0 \\
V6 & 81.08 & 112.5 \\
V7 & 145.97 & 168.28 \\
V8 & 96.48 & 117.5 \\
V9 & 0 & 0 \\
V10 & 0 & 0 \\
V11 & 0 & 0 \\
\hline
\end{tabular}
\end{table}

\begin{table}[!ht]
\centering
\caption{Storage size for 4 representative weeks}
\label{table:storresults_4week}
\begin{tabular}{lll}
\hline
Storage  & 4-week (approach 1) & 4-week (approach 2) \\
\hline 
A & 0 & 0 \\
B & 0 & 0 \\
C & 9.07 & 13.06 \\
D & 86.02 & 100 \\
E & 97.32 & 100 \\
F & 97.32 & 100 \\
G & 0 & 0 \\
H & 0 & 0 \\
I & 100 & 100 \\
K & 100 & 100 \\
L & 0 & 0 \\
M & 100 & 100 \\
N & 100 & 100 \\
O & 0 & 0 \\
P & 100 & 100 \\
Q & 100 & 100 \\
R & 0 & 0 \\
\hline
\end{tabular}
\end{table}

\begin{table}[!ht]
\centering
\caption{Performance for approach 1}
\label{table:rtn4weekres_agg}
\begin{tabular}{llll}
\hline
                 & MADS     & PSO      & Bayesopt \\
\hline
Time (h) & 18.25 & 21.55 & 21.18 \\
$\rho_1$ & 18.9 & 20.34 & 29.89 \\
$\rho_2$ & 26.6 & 30.00 & 27.17 \\
$\rho_3$ & 21.8 & 10.23 & 29.99 \\
$\rho_4$ & 29.3 & 23.98 & 30 \\
$\rho_5$ & 3 & 1.62 & 0.03 \\
$\rho_6$ & 49.6 & 38.73 & 49.96 \\
Best Profit (\$) & -8335.86 & -7987.06 & -8488.58 \\
\hline
\end{tabular}
\end{table}

\begin{table}[!ht]
\centering
\caption{Performance for approach 2}
\label{table:rtn4weekres_dis}
\begin{tabular}{llll}
\hline
                 & MADS     & PSO      & Bayesopt \\
\hline
Time (h) & 19.15 & 22.93 & 31.66 \\
$\rho_1$ & 29.6 & 18.60 & 29.98 \\
$\rho_2$ & 29.9 & 12.02 & 27.96 \\
$\rho_3$ & 20 & 13.69 & 12.33 \\
$\rho_4$ & 30 & 15.36 & 8.21 \\
$\rho_5$ & 9.2 & 0.45 & 4.58 \\
$\rho_6$ & 20 & 2.11 & 0.13 \\
Best Profit (\$) & 7904.33 & 9254.19 & 7645.59 \\
\hline
\end{tabular}
\end{table}

From the results, we find that approach 2 which uses a more detailed high-level model gives a much better solution as compared to approach 1. Furthermore, the time for running the 150 function evaluations from the 3 DFO solvers is lower in approach 1 as compared to that in approach 2. Thus, the results comply with the trend we expect.

\newpage

The full-space model which involves combining the integrated models corresponding to the 4 representative weeks has 52,640 continuous and  14,784 integer variables as well as 67,424 constraints, making it challenging to solve effectively. Therefore solving the integrated with a solver is not a good approach. Additionally, using decomposition methods like the Lagrangian decomposition with each week treated as a subproblem, also fails to converge in a reasonable amount of time because solving even one profile within a day is very difficult. Using metaheuristic methods on the vessel sizes gives a 28-dimensional problem which will take more time as compared to our approach due to the scalability issues in the algorithms. Similar issues will occur with data-driven approaches, which will require a large number of samples due to the large size of the problem. Thus, the PAMSO algorithm proves to be a highly effective method for solving this problem.

\subsection{Integrated planning and scheduling of electrified chemical plants and renewable resources}

In this example, we use our algorithm to solve multi-time scale optimization models for electrification. Electrification, a promising approach for decarbonization, involves powering chemical processes with electricity sourced from renewable green energy rather than fossil fuels. However, these renewable resources have substantial spatial and temporal variations. Additionally, electricity prices fluctuate over time. To fully exploit the economic possibilities of the entire system, it is crucial to account for these fluctuations in the planning of electrified chemical plants and the renewable energy sources that fuel them. The optimal approach involves integrating the planning and scheduling of these plants and resources by formulating a multi-time scale optimization model. 

In this paper, we formulate and solve models for different electric networks. First, we solve a problem involving a network of chemical plants and renewable power-generating units with power lines to connect the locations and the network with external sources of power. Then, we solve a similar problem with the network being isolated.
\subsubsection{Connected network}
In this part, we optimize the planning and scheduling of a low-voltage network of chemical plants and renewable power-generating resources connected to external sources of power.  The multi-time scale optimization model we propose is an MILP model based on the model proposed by \cite{Ramanujam2023DistributedMicrogrid}. The model formulated involves decisions in three time scales: single-time investment decisions, monthly decisions, and hourly decisions. The single-time investment decisions involves the planning decisions to set up the chemical plant network and electric network.  Monthly decisions include the transportation of raw materials and products and inventory management decisions. Representative days are selected for each month. For each representative day, hourly decisions include the power flow in the lines, the operating mode, and the production rate of each plant. The objective of the model is to minimize the total cost including the investment, transportation, inventory, and operating cost. We implement PAMSO on a 5-location case study, the 20-location case study in \Citep{Ramanujam2023DistributedMicrogrid},  and a 200-location case study. While the electric networks in the 5-location and 20-location case study are microgrids, the electric network in the 200-location case study can be treated as a set of connected microgrids.

\textbf{Problem Statement:} We have a few candidate locations where chlor-alkali plants, 1500 kW solar panels and/or 100 kW wind turbines can be placed. We are also given a few locations where the network is connected to an external source of power. Locations in the network can be connected through 12 kW power lines for power to flow between them. The network has a few suppliers of raw materials and consumers of products whose locations are predetermined. Raw materials are obtained from the suppliers and transported to the installed chemical
plants every month. The required chemicals are then produced in the chemical plants and transported every month to consumers whose monthly demand is known. While we try to satisfy the demand of the consumers as much as we can, a penalty is paid for the unsatisfied part of the demand to the consumers. The electrochemical reactions in the plants consume power obtained from the network of renewable power generating units and can also be purchased from the electric utility. In addition, when
the network produces excess electricity, the excess electricity can be sold to the electric utility. Furthermore, we are given the resistance and inductance of the candidate power lines, the electricity price at each hour,
the variable transportation costs, the capital cost of all the generating units, plants, and power lines, and the cost of
different chemicals. We are also given the limits on production rates in different plants, the capacity for power transmission in power lines, and the capacity for power generation of the generating units. We choose 5 representative days for each month. While we constrain the total number of solar panels and wind turbines in the network, we do not specifically constrain the number of chemical plants, solar panels, and wind turbines in each candidate location as well as the number of power lines between 2 locations. Our goal is to optimize the planning decisions (involving the location of the chemical plants and power generating units as well as the connection of power lines), the monthly inventory and transportation decisions, and the hourly operating decisions for the 5 representative days of each month.

\textbf{Strategy:} As mentioned before, we formulate a multi-time scale optimization model to solve this problem. We solve the formulated model using PAMSO. While implementing PAMSO, the high-level model we formulate is an aggregated model based on the original full-space model with hourly operating variables aggregated by month. This high-level model takes decisions in two time scales: single-time investment decisions and monthly decisions (including transportation, inventory, and operating decisions)  and thus has no representative day each month. Therefore, we call it a \emph{0 representative day model} and it outputs the investment decisions. Solving the \emph{0 representative day model} generates the investment decisions, which are then fixed in the full space model to obtain the profit and the operating decisions. The full space model with the investment decisions fixed serves as the low-level model. The low-level model takes decisions in two time scales: monthly transportation and inventory decisions, and hourly operating decisions.  

We parameterize the \emph{0 representative day model} to improve its performance on the full-space model. When we use the \emph{0 representative day model} to generate the investment decisions, we see that solar panels are not profitable in the \emph{0 representative day model}. However, solar panels might be profitable in the full space model. So we add a discount factor to the investment cost of the solar panel \(\rho_1\) as a parameter to the  \emph{0 representative day model}.  Another parameter we add is a discount factor to the capacity of the power lines \(\rho_2\) to account for power losses and AC power flow constraints. Therefore, the parameters which we add to the \emph{0 representative day model} are 
\begin{enumerate}
    \item \(\rho_1\) which represents the discount factor to the investment cost of the solar panel,
    \item \(\rho_2\) which represents the discount factor to the capacity of the power lines.
\end{enumerate}

Both \(\rho_1\) and \(\rho_2\) take values between and including 0 and 1.

\begin{figure}[!ht]
\centering
\includegraphics[width=10cm]{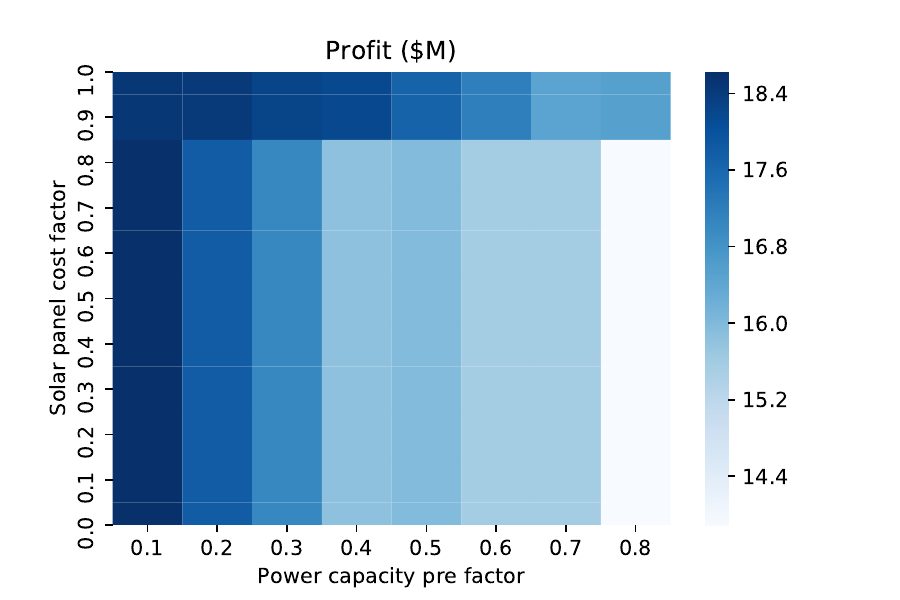}
\caption{Heat map for the MBBF for the 5-location case study}
\label{fig:hm_5loc}
\end{figure}

\textbf{5 location case study}
\\ A case study involving a microgrid with 5 locations is solved. The corresponding full-space model has  313,416 continuous and 86,445 integer variables as well as 1,562,395 constraints. The algorithm is applied to the problem by adding the parameters \(\rho_1\) and \(\rho_2\) to the high-level model.  We have a maximum of 15 solar panels and 150 wind turbines which we can place in the network. To provide an intuitive view of how the selection of the parameters \(\rho_1\) and \(\rho_2\)  affect the multi-scale performance of the high-level model, a heat map is shown in Figure \ref{fig:hm_5loc} with ranges that yield high profit. We notice that the best profit occurs around \(\rho_1 \leq 0.8\) and \(\rho_2 \approx 0.1.\)

During the implementation of PAMSO, we apply the MADS, BayesOpt, and PSO algorithms for 150 function evaluations each to optimize the associated MBBF. While we start with \(\rho_1 = 1\) and \(\rho_2 = 1\) while applying the MADS and PSO algorithms, BayesOpt chooses its own initial parameters. 

\begin{figure}[!ht]
\centering
\includegraphics[width=10cm]{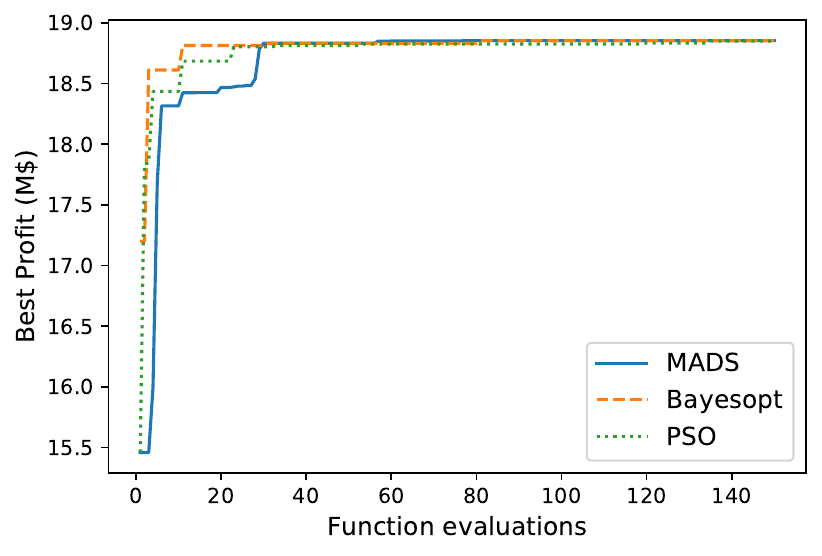}
\caption{Plot for convergence for the 5-location case study}
\label{fig:iter_5loc}
\end{figure}

The convergence of the algorithms for the case study is
shown in Figure \ref{fig:iter_5loc}. From the figure, we can see that MADS reaches its best solution in lesser time as compared to the other 2 DFO solvers. Therefore, we stick with using MADS for the rest of the experiments.

\begin{figure}[H]
\captionsetup[subfigure]{justification=centering}
    \centering
    \begin{subfigure}[b]{1\textwidth}
        \centering
    \includegraphics[width=1\linewidth]{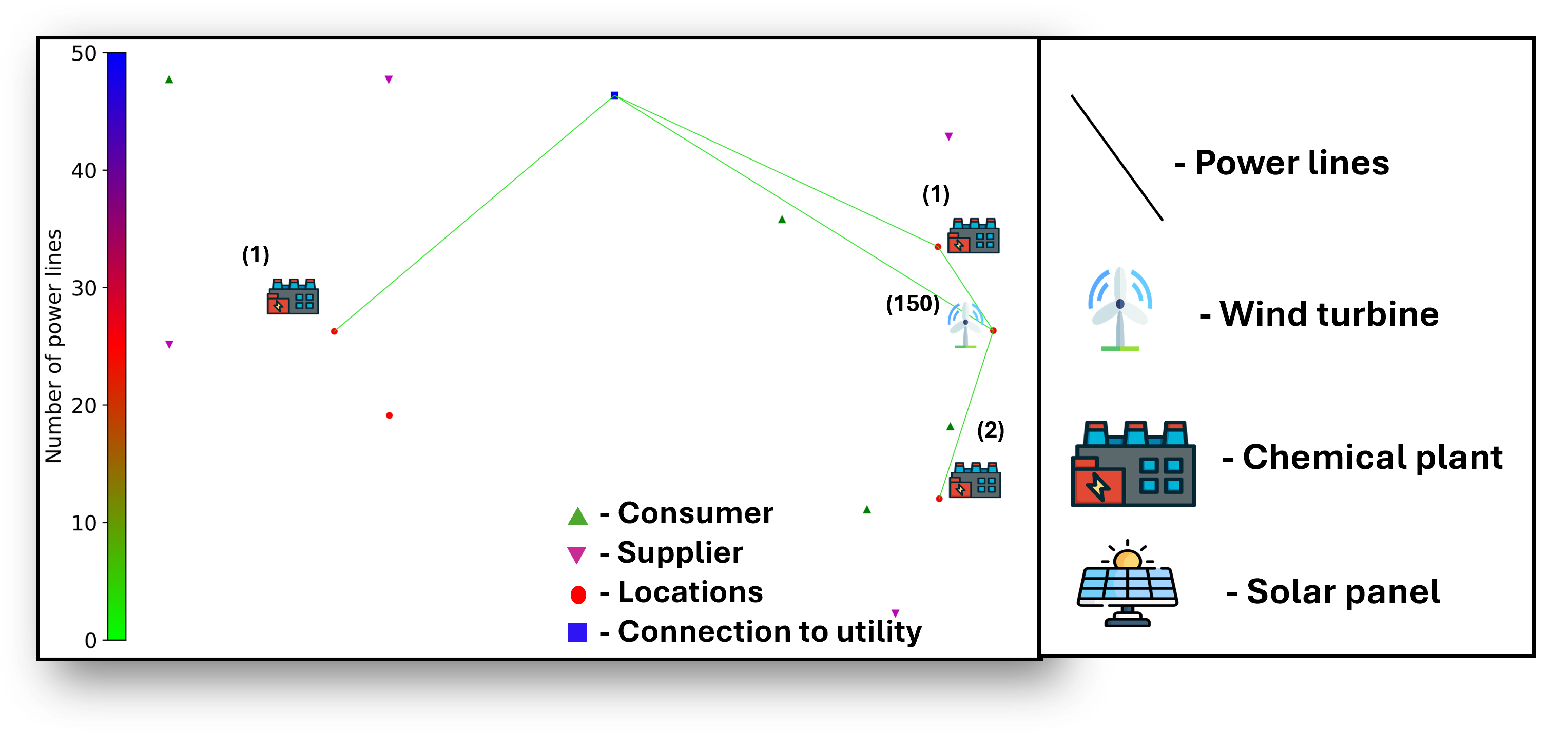}
\caption{Investment decisions from base parameters }

    \end{subfigure}%
\vfill
\vspace{0.25cm}
\vfill
    \begin{subfigure}[b]{1\textwidth}
        \centering
    \includegraphics[width=1\linewidth]{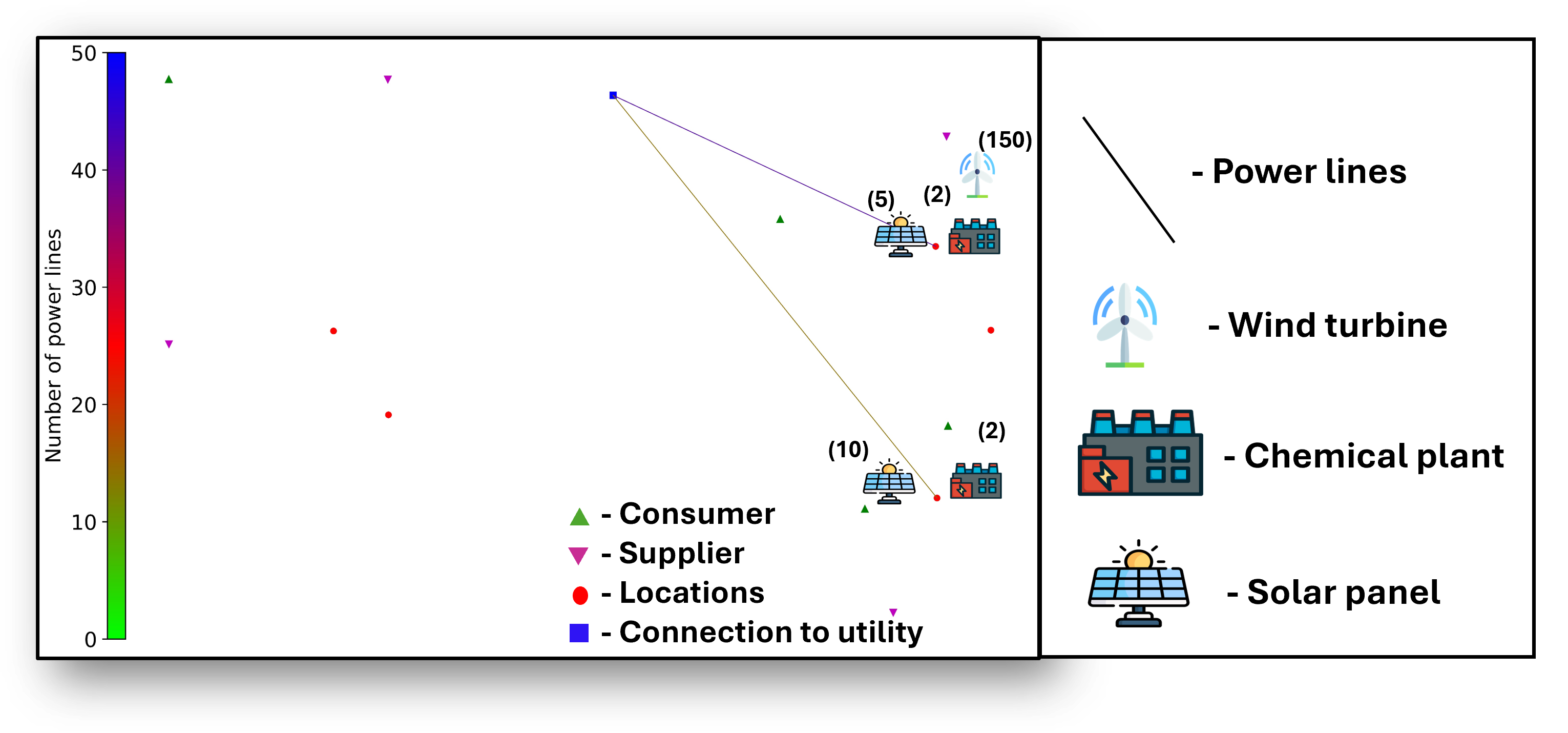}
\caption{Optimal investment decisions from the DFO algorithms}

    \end{subfigure}
     \caption{Investment decisions for the 5-location case study}
 \label{fig:5location}
\end{figure}

The computing time and the optimal parameters are shown in Table \ref{table:rlocres}. The profit at the base parameters is \$15.46 M. The best profit obtained is \$18.85 M with an optimality gap of 3.23\% compared with the LP relaxation (\$19.48M). 
\begin{table}[!ht]
\centering
\caption{Performance of DFO Algorithm}
\label{table:rlocres}
\begin{tabular}{l l l l l}
\hline
DFO Algorithm & Time(h) & \(\rho_1\) & \(\rho_2\) & Best Profit  (M\$)\\
\hline
MADS     & 1.88   & 0.29 & 0.056
  & 18.85  \\
PSO      & 1.85   &0 & 0.055
  & 18.85   \\
Bayesopt & 5.53 & 0.038 & 0.056
 & 18.85 \\

\hline
\end{tabular}
\end{table}
The investment decisions at the base parameters and the best investment decisions at the end of the implementation of the DFO algorithms are shown in Figure \ref{fig:5location}. The main difference between the investment decisions at the base parameters and the optimal investment decisions at the end of the implementation of the 3 DFO solvers is that there are more solar panels and power lines installed in the latter case.

\textbf{20 location case study}
 \\ We use PAMSO to solve the case study involving a microgrid with 20 locations in \cite{Ramanujam2023DistributedMicrogrid}. The corresponding full-space model has  2,094,336 continuous, and 346,080 integer variables as well as 10,111,090 constraints. The algorithm is applied to the problem by adding the parameters \(\rho_1\) and \(\rho_2\) to the high-level model. We apply the MADS algorithm for 150 function evaluations starting with \(\rho_1 = 1\) and \(\rho_2 = 1\)  to optimize the associated MBBF. The profit at the base parameters is \$1.82 M. The best profit obtained is \$18.98 M with an optimality gap of 2.77\% compared to the LP relaxation (\$19.52M). While the time taken for the 150 function evaluations is 17.96 hours, the best solution comes up first in 7.23 hours.

\begin{figure}[H]
\captionsetup[subfigure]{justification=centering}
    \centering
    \begin{subfigure}[b]{0.8\textwidth}
        \centering
    \includegraphics[width=1\linewidth]{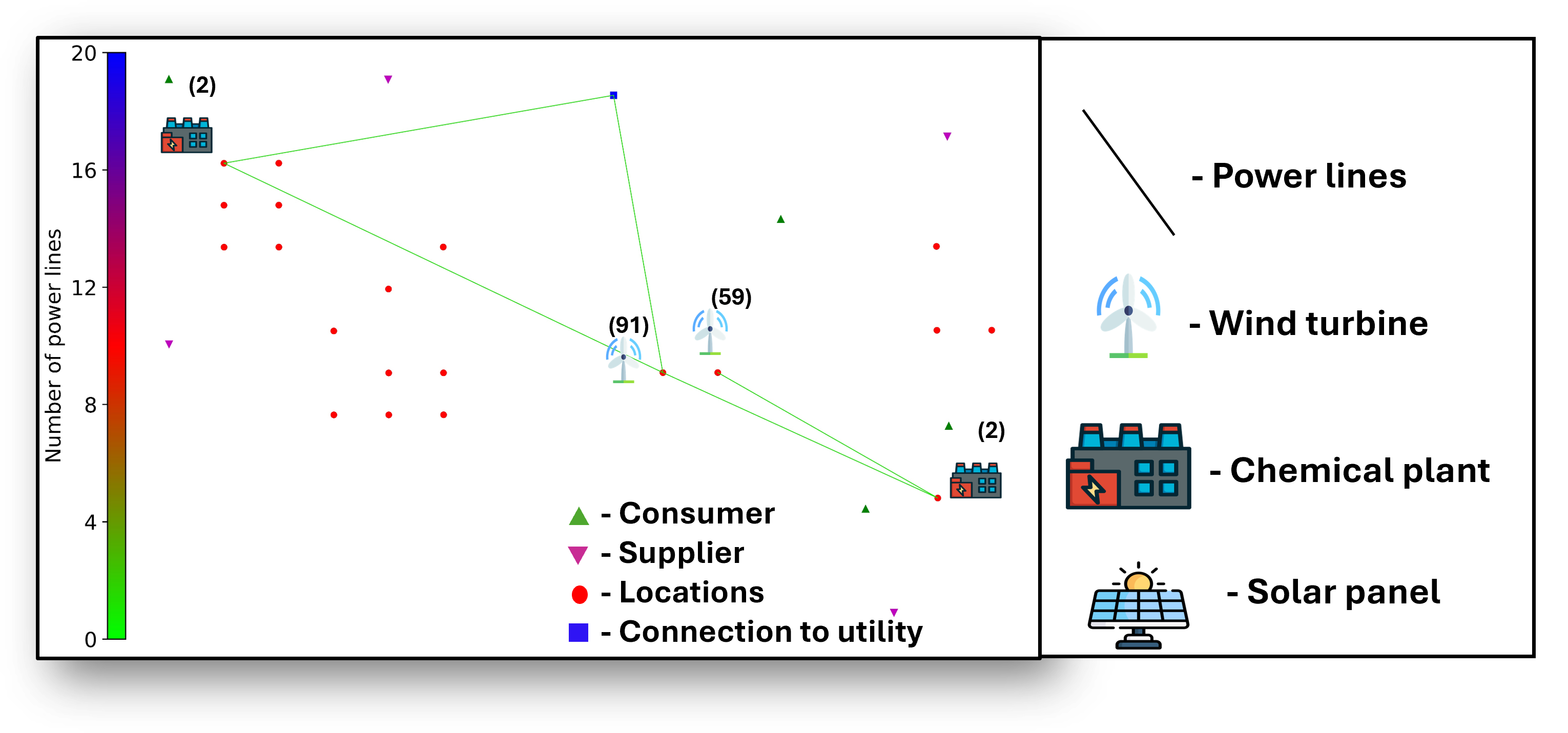}
\caption{Investment decisions from base parameters }

    \end{subfigure}%
\vfill
\vspace{0.25cm}
\vfill
    \begin{subfigure}[b]{0.8\textwidth}
        \centering
    \includegraphics[width=1\linewidth]{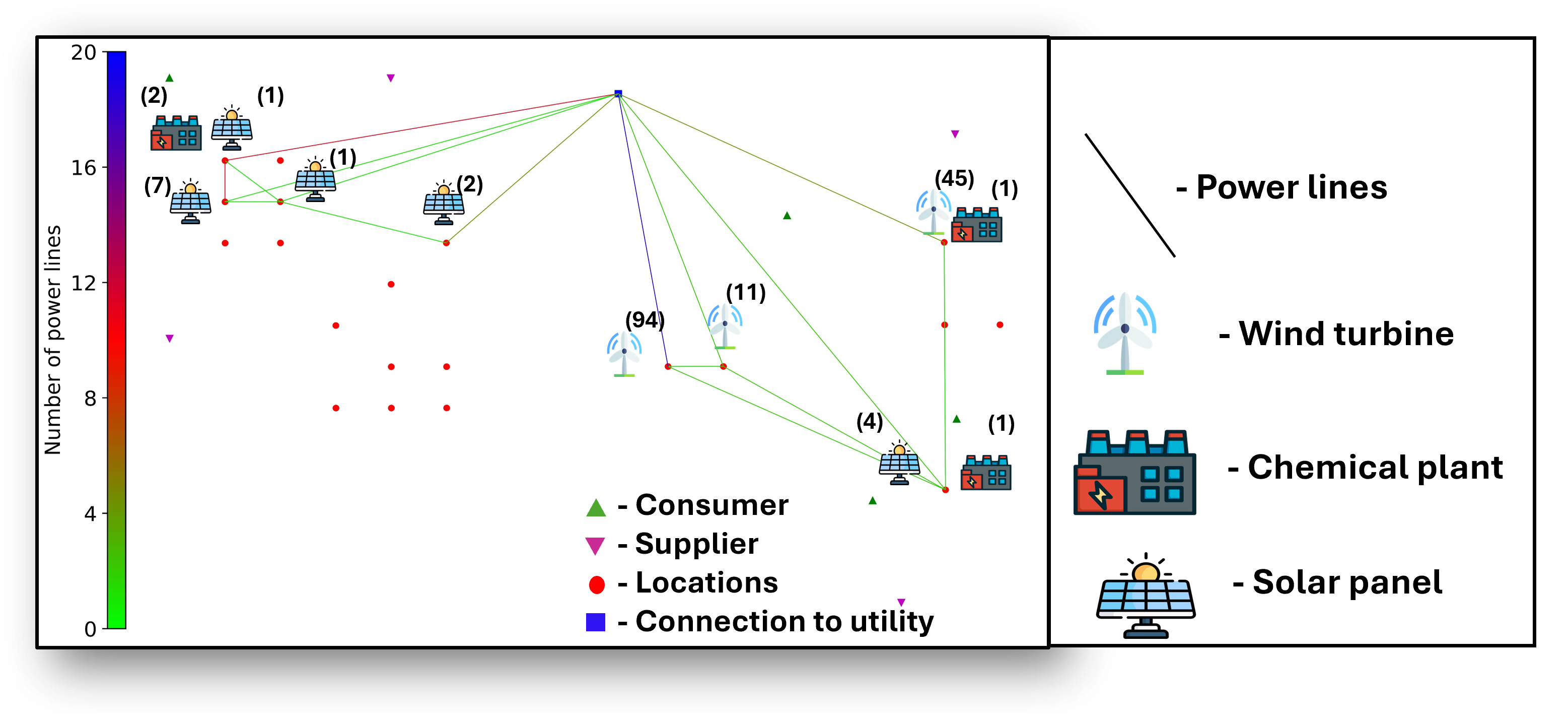}
\caption{Optimal investment decisions from the DFO algorithm}

    \end{subfigure}
     \caption{Investment decisions}
 \label{fig:20location}
\end{figure}

The optimal parameters obtained are \(\rho_1 = 0.1\) and \(\rho_2 = 0.072\). The investment decisions at the base parameters and the end of the implementation of the MADS algorithm are shown in Figure \ref{fig:20location}. The main difference between the investment decisions at the base parameters and the optimal investment decisions at the end of the implementation of the MADS algorithm is that there are more solar panels and power lines installed in the latter case.

\begin{figure}[H]
\captionsetup[subfigure]{justification=centering}
    \centering
    \begin{subfigure}[b]{0.8\textwidth}
        \centering
    \includegraphics[width=1\linewidth]{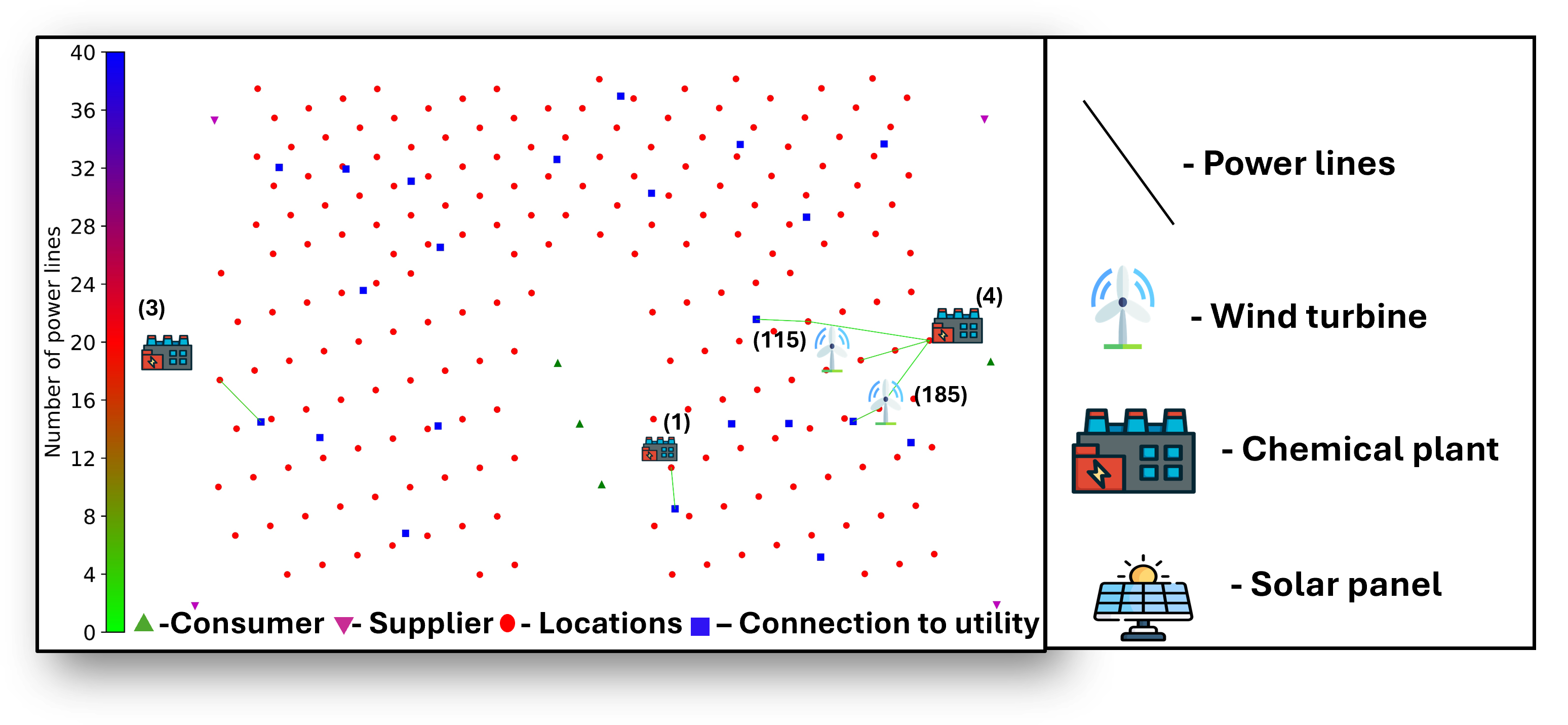}
\caption{Investment decisions from base parameters }

    \end{subfigure}%
\vfill
\vspace{0.25cm}
\vfill
    \begin{subfigure}[b]{0.8\textwidth}
        \centering
    \includegraphics[width=1\linewidth]{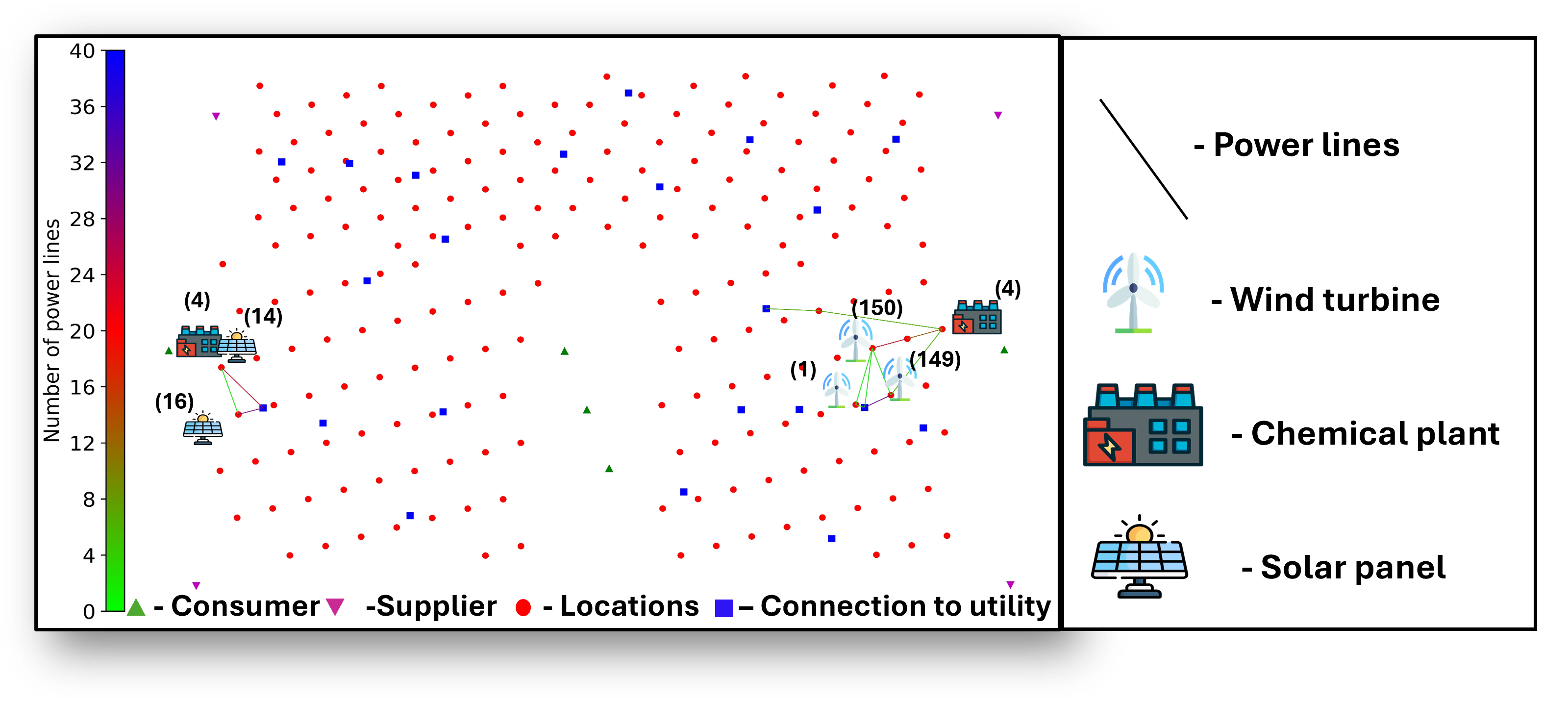}
\caption{Optimal investment decisions from the DFO algorithm}

    \end{subfigure}
     \caption{Investment decisions for the 200-location case study}
 \label{fig:200location}
\end{figure}

\begin{figure}[!ht]
\centering
\includegraphics[width=10cm]{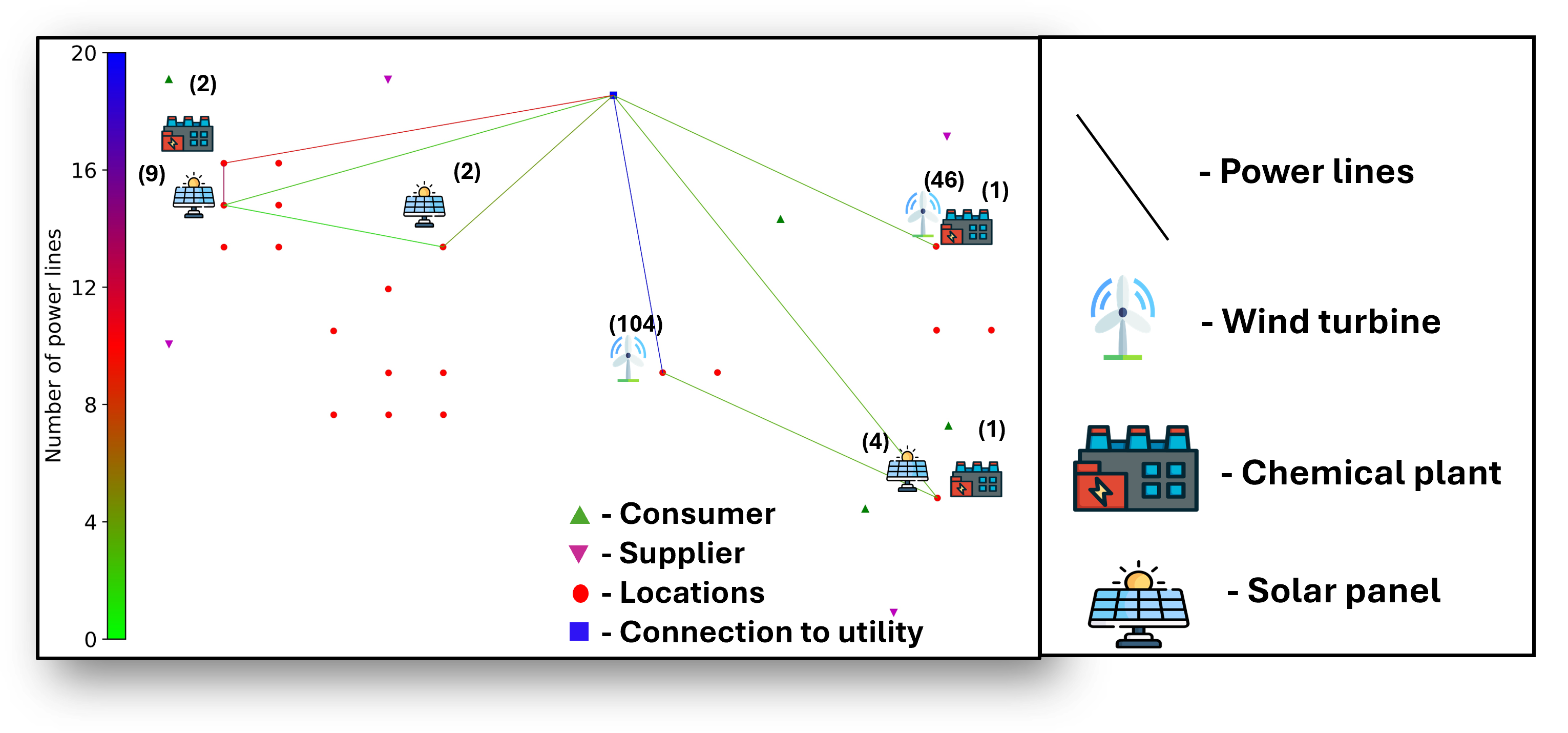}
\caption{Optimal investment decisions obtained from transfer learning of parameters from 5-location case study to 20-location case study}
\label{fig:20_trans}
\end{figure}

 \textbf{200 location case study}
\\ A case study involving a network with 200 locations is solved. The network can be treated as a set of connected microgrids. To take into account the low-voltage characteristics of the network, we connect locations with power lines only if the distance between the two locations is within 35 km. The corresponding full-space model has  22,693,608 continuous and 3,461,396 integer variables as well as 108,920,886 constraints. The algorithm is applied to the problem by adding the parameters \(\rho_1\) and \(\rho_2\) to the high-level model. While implementing the algorithm, we only consider the subset of locations where resources, plants, or power lines are placed while using the low-level model. We ignore the other locations. This significantly reduces the model loading timing and the memory requirements while working with the low-level model. We modify the demand and the network area for the case study. We have a maximum of 30 solar panels and 300 wind turbines which we can place in the network. We apply the MADS algorithm  for 150 function evaluations starting with \(\rho_1 = 1\) and \(\rho_2 = 1\)  to optimize the associated MBBF. The profit at the base parameters is -\$289.91 M. The best profit obtained is \$34.22 M. While the time taken for the 150 function evaluations is 28.02 hours, the best solution comes up first in 24.01 hours. The optimal parameters obtained is \(\rho_1 = 0.84\) and \(\rho_2 = 0.07\). We are unable to run the LP relaxation due to memory issues induced by the large size of the model. The investment decisions at the base parameters and the end of the implementation of the DFO algorithms are shown in Figure \ref{fig:200location}.  The main difference between the investment decisions at the base parameters and the optimal investment decisions at the end of the implementation of the MADS algorithm is that there are more solar panels and power lines installed in the latter case. 
  \textbf{Transfer Learning}
  \\ We apply transfer learning by learning the parameters \(\rho_1\) and \(\rho_2\) from the 5-location example for the 20-location example. To implement transfer learning, we use the parameters obtained from the 5-location case study and apply the MADS algorithm on the associated MBBF for 20 function evaluations in a parameter range of length 0.1, with the initial parameters being the center of the range. The best profit we obtain is \$18.97 M with \(\rho_1 = 0.28\) and \(\rho_2 = 0.08.\) The process takes 2.2 hours. Thus the entire process consisting of 150 function evaluations of the 5-location case study and 20 function evaluations of the 20-location case study gives a solution with a gap of 0.05\% compared to starting from scratch, but it does so in just 56\% of the time required for the initial solution to appear.  The optimum investment decisions are shown in Figure \ref{fig:20_trans}.

Furthermore, we apply transfer learning by learning the parameters \(\rho_1\) and \(\rho_2\) from the 5-location example and 20-location example for the 200-location example in a similar manner. The best profit we obtain from transfer learning from the 5-location case study is \$34.22 M with \(\rho_1 = 0.3\) and \(\rho_2 = 0.07\) and the process takes 3.6 hours. Therefore, the entire process consisting of 150 function evaluations of the 5-location case study and 20 function evaluations of the 200-location case study gives a solution approximately equal to the solution obtained from starting from scratch, but it does so in just 22\% of the time required for the initial solution to appear.

\begin{figure}[H]
\captionsetup[subfigure]{justification=centering}
    \centering
    \begin{subfigure}[b]{0.8\textwidth}
        \centering
    \includegraphics[width=1\linewidth]{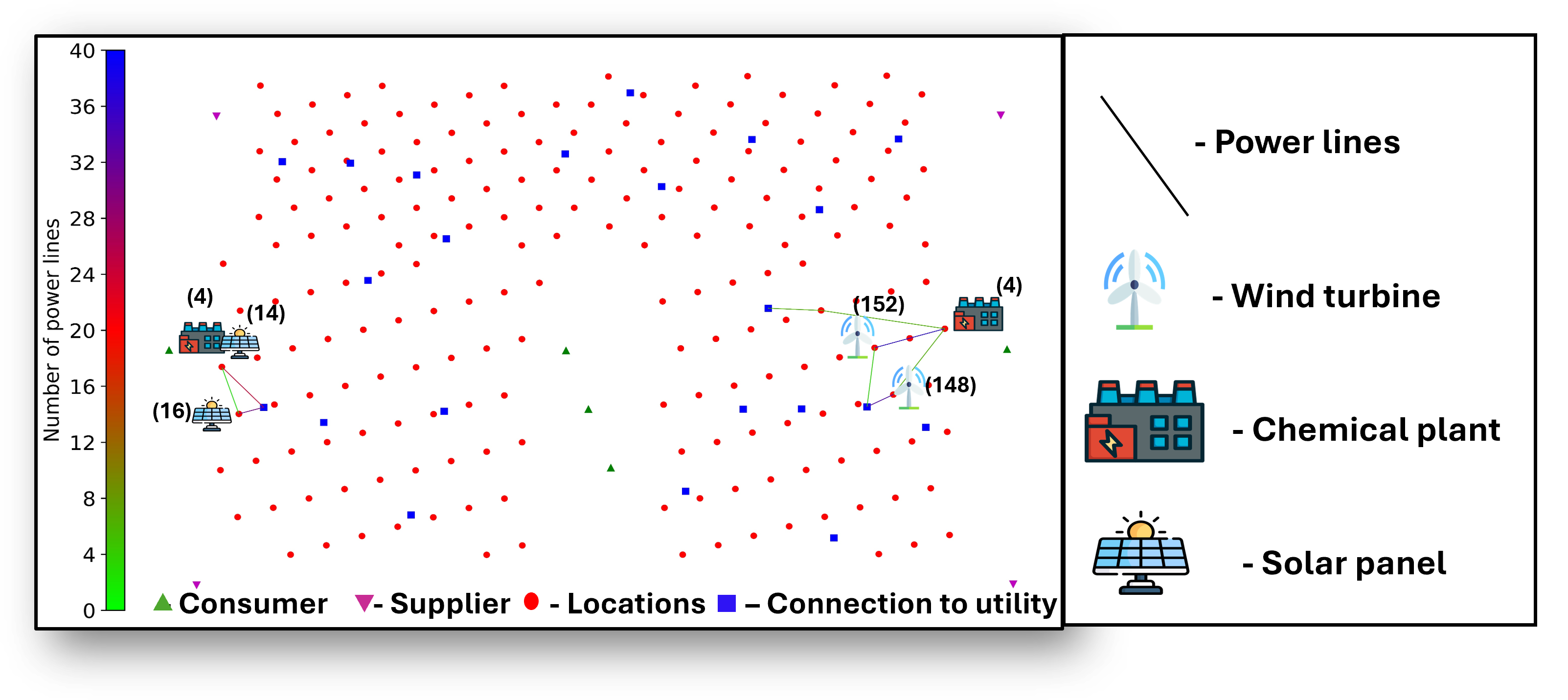}
\caption{Transfer learning of parameters from 5-location case study to 200-location case study }

    \end{subfigure}%
\vfill
\vspace{0.25cm}
\vfill
    \begin{subfigure}[b]{0.8\textwidth}
        \centering
    \includegraphics[width=1\linewidth]{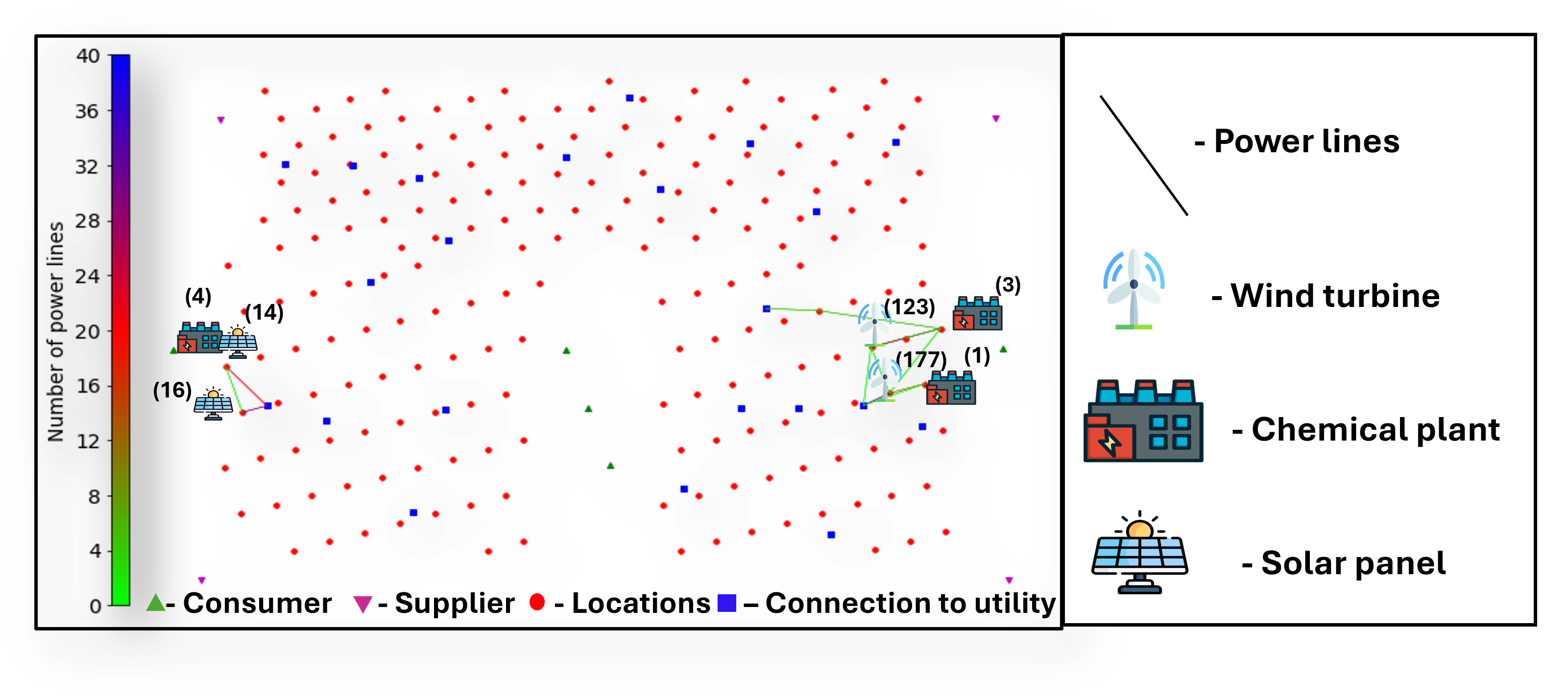}
\caption{Transfer learning of parameters from 20-location case study to 200-location case study}

    \end{subfigure}
     \caption{Optimal investment decisions obtained from transfer learning of parameters to the 200-location case study}
 \label{fig:200_trans}
\end{figure}

The best profit we obtain from transfer learning from the 20-location case study is \$34.15 M with \(\rho_1 = 0.1\) and \(\rho_2 = 0.07\) and the process takes 3.65 hours. The entire process consisting of 150 function evaluations of the 20-location case study and 20 function evaluations of the 200-location case study gives a solution with a gap of 0.2\% compared to starting from scratch, but it does so in 90\% of the time required for the initial solution to appear.  The optimum investment decisions are shown in Figure \ref{fig:200_trans}.

\subsubsection{Isolated microgrid}
In this part, we optimize the planning and scheduling of an isolated microgrid of chemical plants, renewable power-generating units, as well as batteries for electric storage. Instead of using a microgrid connected to the electric utility, we use an isolated microgrid with electric power storage. The case study we solve is similar to the 20-location case study in  \Citep{Ramanujam2023DistributedMicrogrid} with electric storage and no point of common coupling. We use lithium ion batteries of storage capacity 200 kWh for electric storage and obtain the related data from \Citep{Lara2018DeterministicAlgorithm} and \Citep{Schmidt2017TheRates}. We can have a maximum of a thousand  1500kW solar panels and a thousand 100kW wind turbines in the network. The problem can be solved by formulating a multi-time scale optimization model based on the model proposed by \Citep{Ramanujam2023DistributedMicrogrid}.

\textbf{Strategy:}  We solve the problem using the PAMSO algorithm. We modify the \emph{0 representative day model} proposed in the previous part by ensuring that no power is transferred from/to external sources and use the modified model as the high-level model. Solving the high-level model generates the locations of plants, power-generating units, and power lines. These decisions are then fixed in the full space model which acts as the low-level model and provides the profit as well as the monthly transportation and inventory decisions, the hourly operating decisions, and the allocation of the lithium-ion batteries. The lithium-ion batteries are allocated in the low-level model and not in the high-level model, as the high-level model does not accurately capture storage requirements due to the nature of aggregation. The high-level model aggregates variables by month, while electric storage is charged and discharged to a particular value for each day. This implies that, on a monthly basis, the net power obtained from storage after charging and discharging is 0 implying that there is no role for storage in the model. The variables associated with electric storage are integer variables and are easy to optimize after determining the renewable generating units and power lines. Thus we include the allocation of electric storage in the low-level model. 

\begin{figure}[!ht]
\centering
\includegraphics[width=10cm]{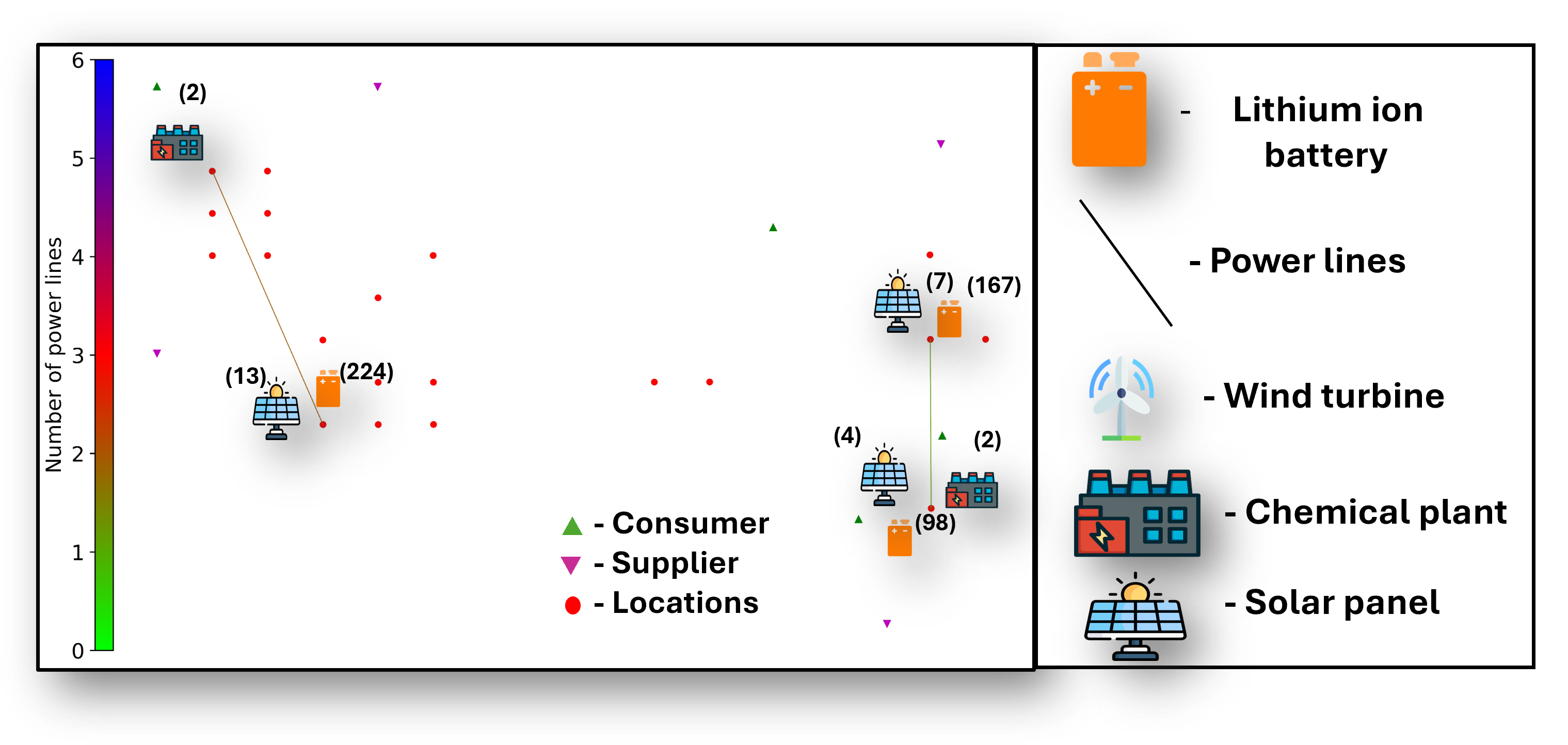}
\caption{Optimal investment decisions without parameters}
\label{fig:20_iso_base_base}
\end{figure}

The temporal variations provided by solar panels and wind turbines are different. Solar panels generate negligible power at the beginning and end of the day, with a peak output around midday. In contrast, wind turbines produce power based on wind speed, which can vary throughout the day. This implies that there are times when solar panels generate more power than wind turbines and vice versa.
While this temporal variation is captured in the full-space model, it is not captured in the high-level model. Without any parameters, the investment decisions obtained from the high-level model are 24 solar panels and no wind turbines. These investment decisions are shown in Figure \ref{fig:20_iso_base_base}. Due to the temporal variations in the power obtained from solar panels, a lot more batteries are required to power the network. The net result includes more investment costs as well as having more of the demand unsatisfied. Wind turbines can provide power in parts of the day when solar panels are not effective. Therefore we would need a mixture of solar panels and wind turbines as well as batteries to power the network well. To ensure that we install wind turbines, we add a minimum number of wind turbines installed  \(\rho_1\) as a parameter to the  \emph{0 representative day model}. The number of solar panels will then be determined by the number of wind turbines and the power demand. Another parameter we add is a discount factor to the capacity of the power lines \(\rho_2\) to account for power losses and AC power flow constraints. Therefore, the parameters which we add to the \emph{0 representative day model} are 
\begin{enumerate}
    \item \(\rho_1\) which represents the minimum number of wind turbines installed in the network
    \item \(\rho_2\) which represents the discount factor to the capacity of the power lines.
\end{enumerate}

While \(\rho_1\) is an integer that takes values between 0 and the maximum possible number of wind turbines in the network, \(\rho_2\) takes values between and including 0 and 1.

\begin{figure}[!ht]
\centering
\includegraphics[width=10cm]{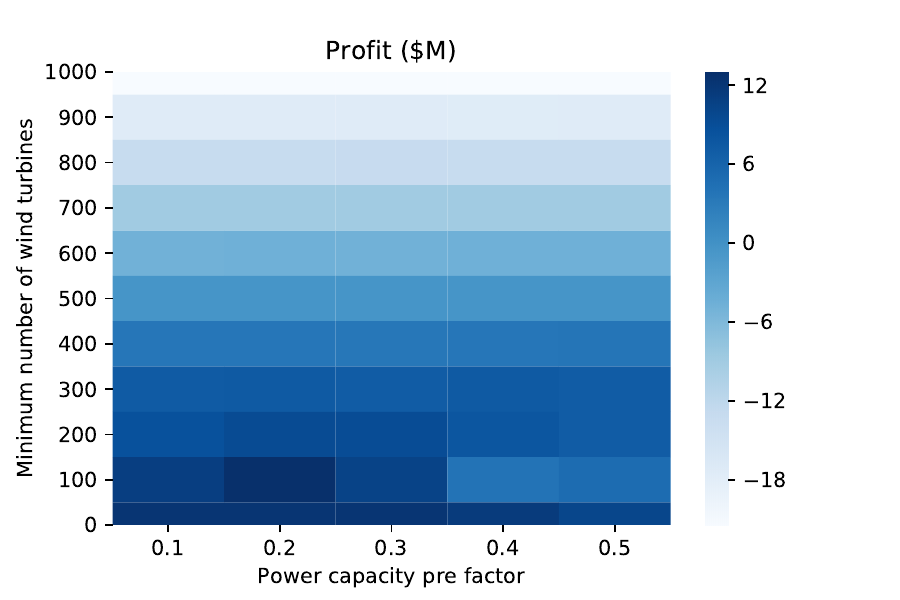}
\caption{Heat map for the MBBF for the isolated case study}
\label{fig:hm_iso}
\end{figure}

\textbf{Results:}
\\ We solve the case study. The corresponding full-space model has  2,180,736 continuous and 346,100 integer integer variables as well as 10,344,110 constraints. The algorithm is applied to the problem by adding the parameters \(\rho_1\) and \(\rho_2\) to the high-level model. To understand the effect of the selection of the parameters \(\rho_1\) and \(\rho_2\) on the multi-scale performance of the high-level model, a heat map is shown in Figure \ref{fig:hm_iso} with ranges that yield high profit. We notice that the best profits occur around the \(\rho_1 \leq 100\) and \(0 \leq \rho_2 \leq 0.2.\)  

\begin{figure}[H]
\captionsetup[subfigure]{justification=centering}
    \centering
    \begin{subfigure}[b]{0.8\textwidth}
        \centering
    \includegraphics[width=1\linewidth]{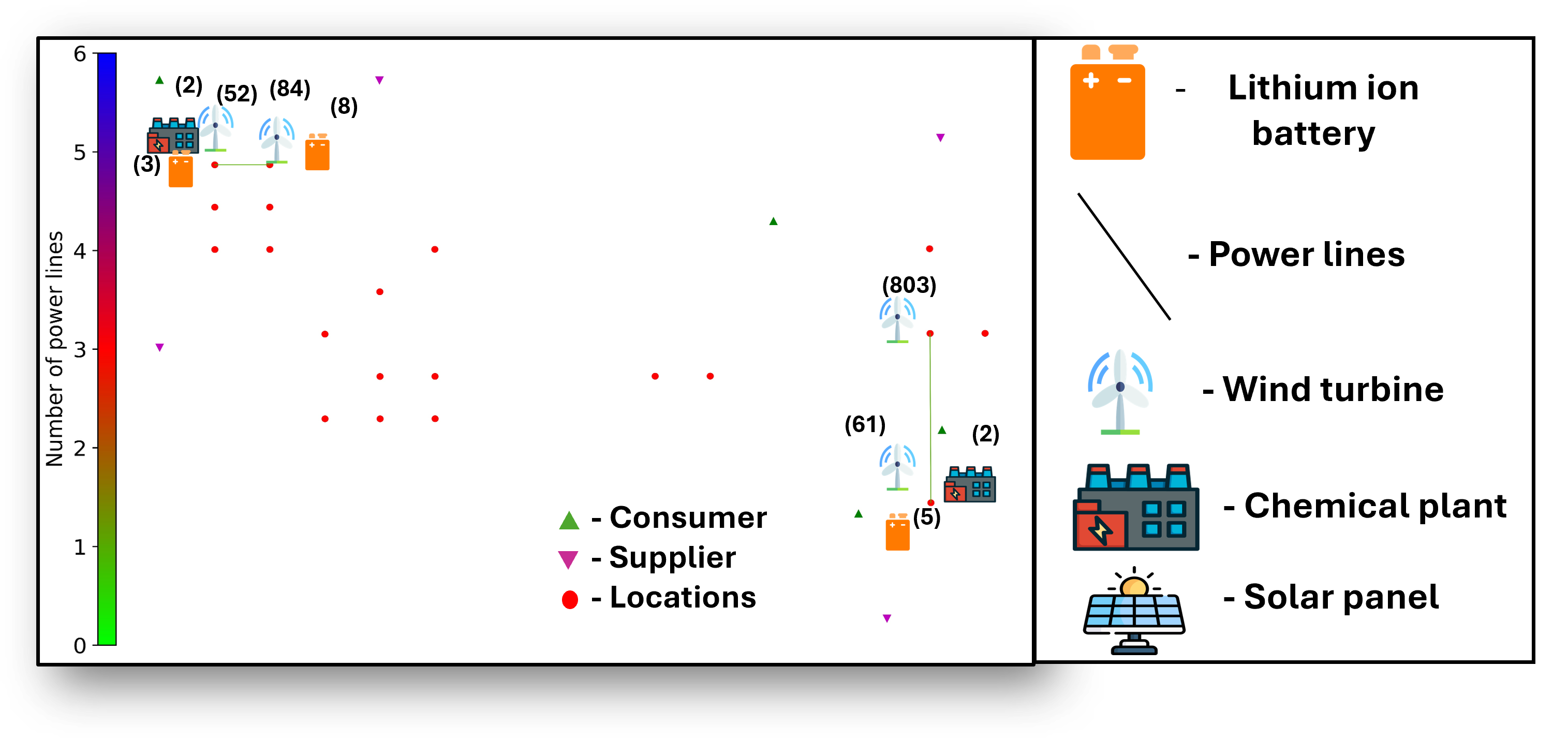}
\caption{Investment decisions from base parameters }

    \end{subfigure}%
\vfill
\vspace{0.25cm}
\vfill
    \begin{subfigure}[b]{0.8\textwidth}
        \centering
    \includegraphics[width=1\linewidth]{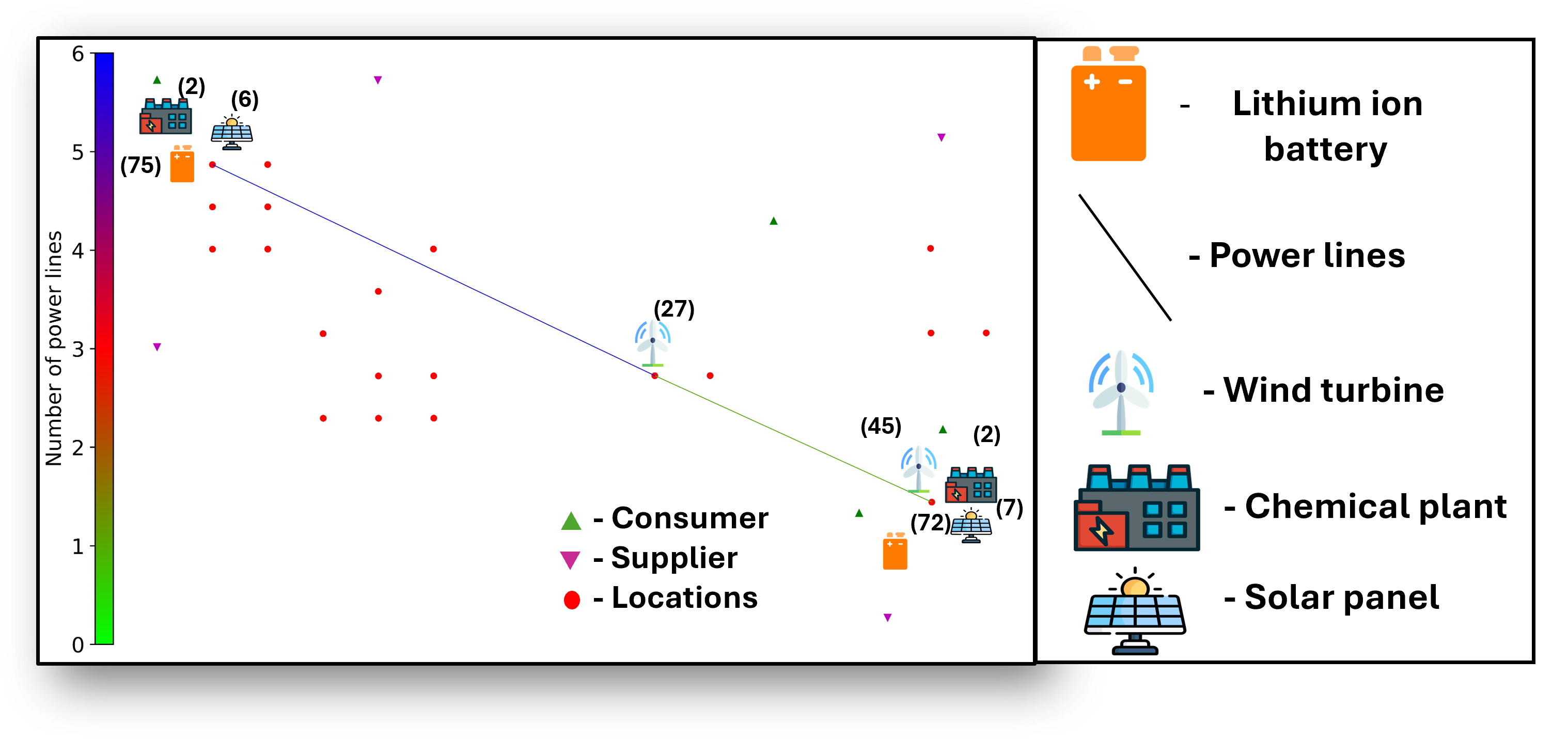}
\caption{Optimal investment decisions from the DFO algorithms}

    \end{subfigure}
     \caption{Investment decisions for the case study with isolated microgrid}
 \label{fig:iso20location}
\end{figure}

We apply the MADS algorithm for 200 function evaluations starting with \(\rho_1 = 1000\) and \(\rho_2 = 1\)  to optimize the associated MBBF. The profit at the base parameters is -\$21.78 M. The best profit obtained is \$14.12 M with an optimality gap of 3.88\% compared to the LP relaxation (\$14.69M). While the time taken for the 200 function evaluations is 34.95 hours, the best solution comes up first in 27.53 hours. The optimal parameters obtained are \(\rho_1 = 72\) and \(\rho_2 = 0.0651\). The investment decisions at the base parameters and at the end of the implementation of the DFO algorithm are shown in Figure \ref{fig:iso20location}. The main difference between the investment decisions at the base parameters and at the end of the implementation of the DFO algorithm is the number of solar panels and wind turbines in the network. 

Isolated networks are more applicable for a moderate number of candidate locations (like 20 locations) close to each other as opposed to a large case study (like 200 locations). This is because isolated networks are powered only by renewable resources. Renewable resources vary in their output across different time scales. There are times of the day when these resources do not provide much power. While having electric storage can help handle these variations to an extent, it is extremely costly to rely only on electric storage for larger networks. Larger networks might also need a large number of power lines to transfer power across large distances making it more costly.   
\subsubsection{Applicability of other algorithms}
For these problems, the Benders decomposition algorithm cannot be applied to this problem as the inventory of consecutive months are coupled, making the problem not decomposable by month or by representative days. Furthermore, the Lagrangian decomposition algorithm is applicable, but takes a large number of iterations to converge and the subproblems are difficult to solve with the growing size of the problem \Citep{Ramanujam2023DistributedMicrogrid}. While the 5-location case study is solvable directly without any decomposition algorithms with some approximations within 24 hours, the 20-location and 200-location case studies are not easily solvable even with approximations within 24 hours. \cite{Ramanujam2023DistributedMicrogrid} proposed an aggregation-disaggregation matheuristic to solve the case studies with a few approximations. While the matheuristic takes 7.5 hours to give a solution of \$18.94 M for the 20-location connected microgrid case study, using PAMSO gives a little better solution of \$18.98 M with the best solution obtained in 7.23 hours.  Using transfer learning from the 5-location problem gives a solution of \$18.97 M in 4.08 hours, thus outperforming the matheuristic's performance. The matheuristic is not scalable to large problems like the 200-location problem. Using metaheuristic methods on the investment decisions gives a problem with a large number of dimensions which will take more time as compared to our approach due to the scalability issues in the algorithms. Similar issues will occur with data-driven approaches, which will require a large number of samples due to the large size of the problem.

\section{Scalability of PAMSO} \label{sec:scalability}
An important advantage of the PAMSO algorithm is its scalability in practice. This is ensured while modelling the high-level and low-level models. The high-level model in practice is a low-fidelity model with not much physics-based complexity and can be solved efficiently. The high-level decisions when fixed in the low-level model, reduce the complexity of the low-level model while using the presolve capability of MILP/MINLP solvers.  In contrast, while using existing decomposition algorithms, we encounter subproblems of large size which are of high fidelity and complexity making them not scalable. Furthermore, metaheuristic methods and data-driven methods also do not scale well to large-size problems as they would need a large number of samples for high-dimensional problems. 

The scalability of PAMSO is further enhanced while using tunable parameters independent of the size of the model. This gives the opportunity of transfer learning parameters from a smaller size problem to a larger size problem which significantly reduces the time to solve the larger size problem. 
\section{Conclusions} \label{sec:conclusion}
In this paper, we proposed the Parametric Autotuning Multi-time Scale Optimization algorithm (PAMSO) to solve multi-time scale optimization models. The algorithm was inspired by parametric cost function approximation (CFA) and the tuning of a controller. In our algorithm, we split the multi-time scale optimization model into a high-level and low-level model and we use parameters in the high-level model to reflect the variations and physics in the low-level model. Transfer learning can be applied to the algorithm by learning the parameters from similar optimization models of smaller sizes by choosing parameters independent of the size of the problem. The algorithm was tested on a set of case studies on the integrated design and scheduling in a resource task network and the integrated planning and scheduling of electrified chemical plants and renewable resources. PAMSO was found to be an efficient option for all the case studies. Transfer-learning of parameters helped reduce the time to solve the larger problems (20-location and 200-location) involving connected networks. 

In the case studies we examined, the high-level decisions were neither temporally interdependent nor variable over time. In other words, the high-level decisions: design decisions in the integrated design and scheduling in a resource task network and investment decisions in the integrated planning and scheduling problem are single-time decisions, that do not have any time-based index associated with them and thus do not vary with any time scale. If the decisions fixed from the high-level model are temporarily interdependent, the likelihood of encountering infeasibilities increases. A large number of parameters may be required to deal with the mismatch in these cases. In the future, we would like to find ways to handle such cases.

 \bibliography{references.bib}





\newpage
\begin{appendices}
\section{Integrated planning and scheduling in RTN}
\label{rtnformulation}
In this section, we show the multi-time scale optimization model for integrated planning and scheduling in an RTN for a week.
\subsection{Notations}
\subsubsection{Indices and Sets}
\begin{description}
    \item  [$i \in \set{I}$]  set of tasks
    \item  [$i \in \set{I}_{r}$]  set of tasks associated with resource \(r \in \set{R}\)
    \item [$r \in \set{R}$] set of resources
    \item  [$r \in \set{R}_{mat}$]  set of materials
    \item  [$r \in \set{R}_{mat}^p$]  set of product materials
     \item  [$r \in \set{R}_{mat}^f$]  set of feed materials
      \item  [$r \in \set{R}_{mat}^i$]  set of intermediate materials
     \item  [$r \in \set{R}_{u}$]  set of vessels
     \item  [$r \in \set{R}_{u}^{f}$]  set of feed task vessels
     \item  [$r \in \set{R}_{u}^{p}$]  set of product task vessels
     \item  [$r \in \set{R}_{u}^{i}$]  set of intermediate task vessels
    \item [$t \in \set{T}_h$] set of hours
    \item [$n \in \set{T}_d$] set of days
    \item  [$i \in \set{U}_{r}$]  set of tasks that can be performed in vessel \(r \in \set{R}_{u}\)
\end{description}
\subsubsection{Discrete variables}
\begin{description}
    \item [$N_{i,t}$] Binary variable to indicate if task \(i \in \set{I} \) starts at time \(t \in \set{T}_h \)
    \item [$N^{agg}_{i,n}$] Integer variable to indicate number of times task \(i \in \set{I} \) starts on day \(n \in \set{T}_d \)
\end{description}
\subsubsection{Continuous variables}
\begin{description}
    \item [$E_{i,t}$] task batch size for task \(i \in \set{I} \) started at time \(t \in \set{T}_h \) (non-negative)
    \item [$E^{agg}_{i,n}$] sum of the batch sizes for task \(i \in \set{I} \) started on day \(n \in \set{T}_d \) (non-negative)
    \item [$\pi_{r,t}$] Flow of resource  \(r \in \set{R} \) in (negative of flow out) to the network at time \(t \in \set{T}_h \) 
    \item [$\pi^{agg}_{r,n}$] Flow of resource  \(r \in \set{R} \) in (negative of flow out) to the network on day \(n \in \set{T}_d \)
    \item [$\phi$]  Total cost
    \item [$\phi_{materials}$]  Total cost from materials
    \item [$\phi_{inv}$]  Total amortized investment cost
    \item [$sl_{r,n}$]  Amount of unfulfilled demand for product \(r \in \set{R}_{mat}^p\) on day \(n \in \set{T}_d \) (non-negative)
     \item [$V^{max}_{r}$] Maximum capacity of vessel \(r \in \set{R}_{u}\) (non-negative)
     \item [$X_{r,t}$] Inventory/Availability for resource \(r \in \set{R}\) at the end of time \(t \in \set{T}_h \) (non-negative)
     \item [$X^{agg}_{r,t}$] Inventory for material \(r \in \set{R}_{mat}\) at the end of day \(t \in \set{T}_d \) (non-negative)
     \item [$X^{max}_{r}$] Maximum capacity of storage for material \(r \in \set{R}_{mat}\) (non-negative)
\end{description}
\subsubsection{Parameters}
\begin{description}
\item[$C^{unit}_{r}$] Ammortized investment cost of vessel for  \(r \in \set{R}_{u}\) with 1 unit capacity for a week
    \item[$C^{stor}_{r}$] Ammortized investment cost of storage for  \(r \in \set{R}_{mat}\) with 1 unit capacity for a week
    \item [$C^{mat}_{r}$] price of material  \(r \in \set{R}_{mat}\) per unit
    \item[$C^{startup}_{i}$] cost of starting task \(i\) once
    \item[$D_{r,n}$] Demand for product \(r \in \set{R}_{mat}^p\) on day \(n \in \set{T}_d \)
    \item [$pen$] Penalty given to consumers for not satisfying their demand
    \item [$\mu_{i,r,\theta}$] constant interaction parameter between task \(i \in \set{I}\) and resource \(r \in \set{R}\) at  time periods \(\theta\) from the start of the task
    \item [$\nu_{i,r,\theta}$] variable interaction parameter between task \(i \in \set{I}\) and resource \(r \in \set{R}\) at  time periods \(\theta\) from the start of the task
   \item [$\nu_{i,r}^{overall}$] overall variable interaction parameter between task \(i \in \set{I}\) and material \(r \in \set{R}_{mat}\)  
    \item [$\tau_i$] time taken for task \(i \in \set{I}\) to be completed
    \item [$v^{min}_i$] minimum fraction of batch capacity for the task batch size for task \(i \in \set{I}\)
     \item [$V^{max,val}_r$] maximum possible value of \(V^{max}_r\)
     \item [$X^{max,val}_r$] maximum possible value of \(X^{max}_r\)

\end{description}
\subsection{Full-space model}
\subsubsection{Resource balance constraint}
The resource balance constraint for each resource \(r \in \set{R}\) is written as follows:
\begin{flalign}
   & X_{r,t} = X_{r,t-1}+\sum_{i \in \set{I}_{r}}\sum_{\theta=0}^{\tau_i} \big(  \mu_{i,r,\theta} N_{i,t-\theta}+\nu_{i,r,\theta} E_{i,t-\theta} \big)+\pi_{r,t} \quad \forall t \geq 1, t \in \set{T}_h, r \in \set{R} 
    \label{Resbalancefp}
\end{flalign}

\subsubsection{Demand satisfaction}
We have the mass balance constraints with respect to the demand $D_{r,n}$ and its satisfaction below. The part of the demand not met for the consumer, \(sl_{r,n}\) is penalized in the objective function. 

\begin{flalign}
   & -\left(\sum_{t=24(n-1)+1}^{24n} \pi_{r,t} \right)= D_{r,n} - sl_{r,n} \quad \forall n \geq 1, n\in \set{T}_d,r \in \set{R}^{p}_{mat}
    \label{Dembalancefp}
\end{flalign}

\subsubsection{Resource limit constraint}
The minimum inventory of resources is 0 with the following constraint written:
\begin{flalign}
   & X_{r,t} \geq 0 \quad  \forall t \geq 1, t \in \set{T}_h, r \in \set{R} 
    \label{Resfpmin}
\end{flalign}

The maximum availability for each vessel is 1 and the maximum inventory for each material is \(X^{max}_{r}.\) The following constraints are written to ensure this.

\begin{flalign}
   & X_{r,t} \leq 1 \quad  \forall t \geq 1, t \in \set{T}_h, r \in \set{R}_u 
    \label{Resfpmax1}
\end{flalign}

\begin{flalign}
   & X_{r,t} \leq X^{max}_{r} \quad  \forall t \geq 1, t \in \set{T}_h, r \in \set{R}_{mat} 
    \label{Resfpmax2}
\end{flalign}

\subsubsection{Batch constraint}
The maximum task batch size for task \(i \in \set{U}_{r}\) depends on the capacity of the vessel \(V^{max}_{r}\) associated with it and is constrained in the following way:
\begin{flalign}
   & E_{i,t} \leq V^{max}_{r}N_{i,t} \quad  \forall t \geq 1, t \in \set{T}_h, i \in \set{U}_r,r \in \set{R}_{u} 
    \label{batfpmax1}
\end{flalign}

The minimum task batch size for task \(i \in \set{U}_{r}\) would depend on the capacity of the vessel \(V^{max}_{r}\)  associated with it as well as \(v^{min}_i\) and is constrained in the following way:

\begin{flalign}
   & E_{i,t} \geq v^{min}_iV^{max}_{r}N_{i,t} \quad  \forall t \geq 1, t \in \set{T}_h, i \in \set{U}_r,r \in \set{R}_{u} 
    \label{batfpmin1}
\end{flalign}

\subsubsection{Additional constraints}
Additional constraints are added to make the model more meaningful and easier to solve. We add constraints on \(\pi\) based on the type of resource we are dealing with i.e., products always move out of the system, feed moves into the system, and intermediates as well as vessels do not move in or out of the system

\begin{flalign}
   & \pi_{r,t} \leq 0 \quad  \forall t \geq 1, t \in \set{T}_h, r \in \set{R}_{mat}^{p} 
    \label{add1}
\end{flalign}

\begin{flalign}
   & \pi_{r,t} \geq 0 \quad  \forall t \geq 1, t \in \set{T}_h, r \in \set{R}_{mat}^{f} 
    \label{add2}
\end{flalign}

\begin{flalign}
   & \pi_{r,t} = 0 \quad  \forall t \geq 1, t \in \set{T}_h, r \in \set{R}_{mat}^{i} \cup  \set{R}_{u}
    \label{add3}
\end{flalign}
The initial amount of materials is 0. The initial availability of each vessel is 1.
\begin{flalign}
   & X_{r,0} = 0 \quad  \forall  r \in \set{R}_{mat}
    \label{add6}
\end{flalign}

\begin{flalign}
   & X_{r,0} = 1 \quad  \forall  r \in \set{R}_{u}
    \label{add7}
\end{flalign}

We also add bounds on \(V^{max}_{r}, X^{max}_{r}\) by estimating \(V^{max,val}_{r},X^{max,val}_{r}\)  from the demand and requirements of the problem.

\begin{flalign}
   & 0 \leq V^{max}_{r} \leq V^{max,val}_{r} \quad  \forall r \in  \set{R}_{u}
    \label{add4}
\end{flalign}

\begin{flalign}
   & 0 \leq X^{max}_{r} \leq X^{max,val}_{r} \quad  \forall  r \in \set{R}_{mat}
    \label{add5}
\end{flalign}
\subsubsection{Objective function}
We maximize the net profit or minimize the net cost of the system.

The amortized investment cost can be written as follows
\begin{flalign}
    & \phi_{inv} = \sum_{r \in R_u}C^{unit}_{r}(V^{max}_{r})^{0.6} +  \sum_{r \in R_{mat}}C^{stor}_{r}(X^{max}_{r})^{0.6} 
    \label{invcost}
\end{flalign}

 We try to meet the demand of each consumer. For the part of the demand unmet $sl_{r,n}$, we pay a penalty of $pen \cdot C^{mat}_{r}sl_{r,n}.$ The net cost of materials can be written as follows

\begin{flalign}
    & \phi_{materials} = \sum_{r \in R_{mat}}\sum_{t \in \set{T}_h}C^{mat}_{r}\pi_{r,t} + \sum_{r \in R_{mat}^p}\sum_{n \in \set{T}_d} pen \cdot C^{mat}_{r}sl_{r,n}
    \label{matcost}
\end{flalign}

We also add a cost of \(C^{startup}_{i}\) for starting of task \(i.\) The total startup cost is \(\sum_{i \in I}\sum_{t \in \set{T}_h}C^{startup}_{i}N_{i,t}.\)

The entire optimization model is summarized in \eqref{eq:fullmodel}
\begin{subequations} \label{eq:fullmodel}
    \begin{align}
        \min \quad & \phi =  \sum_{i \in I}\sum_{t \in \set{T}_h}C^{startup}_{i}N_{i,t}+\phi_{material}+\phi_{inv} & \\
        \text{s.t. } \quad & \eqref{Resbalancefp} - \eqref{matcost},  
    \end{align}
\end{subequations}

\subsection{High-level model}
The high-level model is adapted from the full space model with the scheduling decisions aggregated on a daily basis. That is instead of having variables for each hour, we have variables for a day. With a slight abuse of notations, the variables and their aggregated version are summarized in Table \ref{table:rtnagg}.

\begin{table}[H]
\centering
\caption{Parameters for implementing PAMSO on integrated design and scheduling model}
\label{table:rtnagg}
\begin{tabularx}{\textwidth}{>{\centering\arraybackslash}X>{\centering\arraybackslash}X>{\centering\arraybackslash}X>{\centering\arraybackslash}X}
\hline
Variable in full-space & Physical meaning & Variable in high-level model & {Physical   meaning} \\
\hline \\
$N_{i,t}$ & Binary variable to indicate if task \(i \in \set{I}\) is started in time \(t \in   \set{T}_h \) & $N^{agg}_{i,n}$ & {Integer variable to indicate number   of times task \(i \in \set{I} \) starts on day \(n \in \set{T}_d \)} \\
\hline \\
$E_{i,t}$ & batch size for task \(i \in \set{I} \) started at time \(t \in   \set{T}_h \) (non-negative) & $E^{agg}_{i,n}$ & {sum of the batch sizes for task \(i   \in \set{I} \) started on day \(n \in \set{T}_d \) (non-negative)} \\
\hline \\
$\pi_{r,t}$ & Flow of resource  \(r \in \set{R}   \) in (negative of flow out) to the network at time \(t \in \set{T}_h \) & $\pi^{agg}_{r,n}$ & {Flow of resource  \(r \in   \set{R} \) in (negative of flow out) to the network on day \(n \in \set{T}_d   \)} \\
\hline \\
$X_{r,t}$ & Inventory for resource \(r \in \set{R}\) at the end of time \(t \in   \set{T}_h \) (non-negative) & $X^{agg}_{r,t}$ & {Inventory for material \(r \in   \set{R}_{mat}\) at the end of day \(t \in \set{T}_d \) (non-negative)}  \\
\hline

\end{tabularx}
\end{table}
The overall variable interaction parameter \(\nu^{overall}_{i,r}\) between task \(i \in \set{I}\) and material \(r \in \set{R}_{mat}\)  is the sum of the variable interaction parameter \(\nu_{i,r,\theta}\) over the lifetime of the task \(i\) i.e., \(\nu^{overall}_{i,r} = \sum_{\theta = 0}^{\tau_{i}}\nu_{i,r,\theta}.\)
\newpage
The optimization model is given below:

\begin{subequations} \label{eq:fullmodel1}
    \begin{flalign}
         \min \quad   & \phi = \sum_{i \in I}\sum_{n \in \set{T}_d}C^{startup}_{i}N^{agg}_{i,n} +\phi_{material}+\phi_{inv} \\ 
        \text{s.t. } \quad & X^{agg}_{r,n} = X^{agg}_{r,n-1}+\sum_{i \in \set{I}_{r}} \big(\nu^{overall}_{i,r} E_{i,n} \big)+\pi^{agg}_{r,n}  \quad \forall n \geq 1, n \in \set{T}_d, r \in \set{R}_{mat} \\
        & -\pi^{agg}_{r,n}  = D_{r,n} - sl_{r,n} \quad \forall n \geq 1, n\in \set{T}_d,r \in \set{R}^{p}_{mat} \\  
   & 0 \leq X^{agg}_{r,n} \leq X^{max}_{r} \quad  \forall n \geq 1, n \in \set{T}_d, r \in \set{R}_{mat}  \\
   & v^{min}_{i} V^{max}_{r}N_{i,n} \leq E^{agg}_{i,n} \leq V^{max}_{r}N_{i,n} \quad  \forall t \geq 1, t \in \set{T}_h, i \in \set{U}_r,r \in \set{R}_{u}\\
   & \pi^{agg}_{r,n} \leq 0 \quad  \forall n \geq 1, n \in \set{T}_d, r \in \set{R}_{mat}^{p} \\ 
      & \pi^{agg}_{r,n} \geq 0 \quad  \forall n \geq 1, n \in \set{T}_d, r \in \set{R}_{mat}^{f} \\
     & \pi^{agg}_{r,n} = 0 \quad  \forall n \geq 1, n \in \set{T}_d, r \in \set{R}_{mat}^{i} \cup  \set{R}_{u} \\
   & X_{r,0} = 0 \quad  \forall  r \in \set{R}_{mat} \\
    & X_{r,7} = 0 \quad  \forall  r \in \set{R}_{mat} \\
   & 0 \leq V^{max}_{r} \leq V^{max,val}_{r} \quad  \forall r \in  \set{R}_{u} \\ 
   & 0 \leq X^{max}_{r} \leq X^{max,val}_{r} \quad  \forall  r \in \set{R}_{mat} \\
     & \phi_{inv} = \sum_{r \in R_u}C^{unit}_{r}(V^{max}_{r})^{0.6} +  \sum_{r \in R_{mat}}C^{stor}_{r}(X^{max}_{r})^{0.6} \\
      & \phi_{material} = \sum_{r \in R_{mat}}\sum_{n \in \set{T}_d}C^{mat}_{r}\pi^{agg}_{r,n} + \sum_{r \in R_{mat}^p}\sum_{n \in \set{T}_d} pen \cdot 
 C^{mat}_{r}sl_{r,n}\\
   & \sum_{i \in \set{U}_{r}} N^{agg}_{i,n}\tau_{i} \leq 24 \quad  \forall  n \in \set{T}_d,r \in \set{R}_u 
       \end{flalign}
\end{subequations}

We can write the high-level model defined above in the form defined in the algorithm section as shown below:

\begin{subequations}
    \begin{align}
        \min_{\bm{x},\bm{z}} \quad f(\bm{x})+Q(\bm{x},\bm{z},\bm{\Tilde{\theta}}) & \\
        \text{s.t.} \quad g(\bm{x}) \leq 0 & \\
        H(\bm{x}, \bm{z}, \bm{\Tilde{\theta}}) \leq 0 & 
    \end{align}
\end{subequations}
with \(\bm{x} = [V^{max}_{r},X^{max}_{r}]\) and \(\bm{z}\) being the set of other variables. The tunable parameters are shown in Table \ref{table:rtnparameters}.
 The parametrized cost would be as follows:

 \begin{subequations}
      \begin{flalign}
     & \phi = \sum_{i \in I}\sum_{n \in \set{T}_d} \rho_6 C^{startup}_{i}N^{agg}_{i,n} +\phi_{material}+\phi_{inv} \\
& \phi_{material} = \sum_{r \in R_{mat}}\sum_{n \in \set{T}_d}C^{mat}_{r}\pi^{agg}_{r,n} + \sum_{r \in R_{mat}^p}\sum_{n \in \set{T}_d} \rho_1 pen \cdot C^{mat}_{r}sl_{r,n}
     \end{flalign}
     \begin{equation}
         \begin{split}
              & \phi_{inv} = \sum_{r \in \set{R}_u^f}\rho_2C^{unit}_{r}(V^{max}_{r})^{0.6}+\sum_{r \in \set{R}_u^p}\rho_3C^{unit}_{r}(V^{max}_{r})^{0.6}+ \\ & \sum_{r \in \set{R}_u^i}\rho_4C^{unit}_{r}(V^{max}_{r})^{0.6} +  \sum_{r \in R_{mat}}\rho_5C^{stor}_{r}(X^{max}_{r})^{0.6} \\
         \end{split}
     \end{equation}
 \end{subequations}

\end{appendices}

\end{document}